\documentclass[hidelinks,onefignum,onetabnum]{siamart251216}

\usepackage{amssymb,mathtools,bm}
\usepackage{booktabs,multirow,array}
\usepackage{adjustbox}
\usepackage{graphicx}
\usepackage{algpseudocode}
\usepackage{enumitem}
\usepackage{microtype}

\newsiamthm{assumption}{Assumption}
\newsiamremark{remark}{Remark}
\crefname{assumption}{Assumption}{Assumptions}
\crefname{remark}{Remark}{Remarks}

\newcommand{\kappae}{\kappa}

\newcommand{\spanop}{\operatorname{span}}

\headers{Learned Local Spaces for Multiscale Elements}{Y. Wang, C. X. Kong, and G. Lin}

\title{Multiscale LegONet: Learning Coefficient-Dependent Local Spaces for Composable Multiscale Finite Elements with Localized Error Bounds\thanks{Submitted to the editors \today.
\funding{This work was supported by the DARPA MATHBAC program under grant HR0011-26-3-E152; by the National Science Foundation under grants DMS-2533878, DMS-2053746, DMS-2134209, ECCS-2328241, CBET-2347401 and OAC-2311848; by the U.S.~Department of Energy (DOE) Office of Science Advanced Scientific Computing Research program under the ``Uncertainty Quantification for Multifidelity Operator Learning (MOLUcQ)'' project (Project No.~81739), DE-SC0023161 and the SciDAC LEADS Institute; and by DOE Fusion Energy Sciences under grant DE-SC0024583.}}}

\author{Yueqi Wang\thanks{Department of Mathematics, Purdue University, West Lafayette, IN 47907
  (\email{wang7406@purdue.edu}, \email{kong146@purdue.edu}). These authors contributed equally to this work.}
\and Cindy Xiangrui Kong\footnotemark[2]
\and Guang Lin\thanks{Corresponding author. Department of Mathematics and School of Mechanical Engineering, Purdue University, West Lafayette, IN 47907
  (\email{guanglin@purdue.edu}).}}

\ifpdf
\hypersetup{
  pdftitle={Multiscale LegONet: Learning Coefficient-Dependent Local Spaces for Composable Multiscale Finite Elements with Localized Error Bounds},
  pdfauthor={Y. Wang, C. X. Kong, and G. Lin}
}
\fi

\begin{document}

\maketitle

\begin{abstract}
We introduce multiscale LegONet (Ms-LegONet), a plug-and-play learned block for heterogeneous elliptic problems. One network per spatial dimension is trained once on local coefficient patches and reused, without retraining, across coefficient layouts, sources, boundary data and domain sizes within the trained coefficient class. On a fixed primitive patch the network maps the local coefficient to a zero-boundary space, the only learned object. The multiscale basis functions are the harmonic extensions of a fixed set of boundary traces, each a Galerkin solve in that space. Shared interface traces enforce continuity, larger blocks are assembled by recursive static condensation without any parent-scale network, and a localized C\'ea bound attributes the error of the learned spaces primitive by primitive. The training labels are harmonic responses of single patches, so no global solve enters the training. Against MsFEM, GMsFEM, the neural-operator element method and an end-to-end Fourier neural operator, Ms-LegONet is the most accurate at equal unknowns: its $L^2$ error is $4$ to $150$ times below the baselines on two-dimensional mosaics of $256$ high-contrast primitives, and its median error stays within $0.1$ percentage points of that of exact local spaces in three-dimensional layered media. For each new coefficient it is $1.6$ to $2.3$ times faster than a fine-grid direct solve in two dimensions. In three dimensions it computes the many-query port response of a chain $30$ to $900$ times faster than finite element solves on the same or a finer mesh, at a relative error of $10^{-4}$ against their $6\%$ and $1.6\%$.
\end{abstract}



\section{Introduction}
\label{sec:intro}

Elliptic equations with heterogeneous coefficients govern porous media flow, composite materials and heat conduction \cite{bear2013dynamics,torquato2002random,bensoussan2011asymptotic,jikov2012homogenization}.  Features spread over many scales force a fine discretization, and its cost is paid again for every new coefficient.  This is the bottleneck of many-query settings such as uncertainty quantification and upscaling \cite{ghanem2003stochastic,kuo2016application,benner2015survey}.

Multiscale finite element methods replace fine global trial spaces with basis functions from local coefficient-dependent cell problems \cite{hou1997msfem}.  Spectral enrichment and the generalized multiscale finite element method (GMsFEM) extend the construction to high contrast with modes of local eigenproblems \cite{efendiev2011spectral,efendiev2013gmsfem}, and reduced-basis multiscale methods reuse low-dimensional spaces to accelerate the recurring local solves \cite{abdulle2012rbfehmm,ohlberger2015lrbms,buhr2017arbilomod}.  Component-based static condensation with conforming ports and reduced interior spaces is the setting of the static condensation reduced basis element method \cite{huynh2013scrbe}, with optimal port spaces in \cite{smetana2016optimal}.  In every case the one expensive, coefficient-dependent step is local and repeated for each coefficient; everything after is deterministic and sees the coefficient only through the local spaces.  That local step is the right target for a learned surrogate, which is most reusable when it is independent of the global layout and composable without retraining.

Neural operators learn maps between function spaces; DeepONet, the Fourier neural operator (FNO), the neural-operator framework, physics-informed operator training and multigrid-inspired parameterizations all target the global solution operator \cite{lu2021deeponet,li2021fno,kovachki2023neuraloperator,li2024pino,he2024mgno}.  On domains assembled from independent high-contrast patches that global map has little low-frequency structure to exploit: an end-to-end FNO trained on mosaics of the same patches as our local block stalls at $38\%$ $L^2$ error in two dimensions (\Cref{subsec:2d-diffusion}).
A second line of work learns local or unresolved components, from coarse closures to blocks assembled by domain decomposition, and keeps the known structure deterministic \cite{lutjens2022mno,li2024lno,fang2024boundary,chen2025locality,huang2026sni,he2026selfcompose}.  Three methods reuse a local neural object inside a variational computation.  NOEM embeds a neural operator as a subdomain element whose input includes the boundary state \cite{ouyang2026noem}; GMsFEM-NO replaces the local eigenproblems of GMsFEM by inference of the spectral subspace \cite{rudikov2025gmsfemno}; our earlier LegONet composes independently pretrained, plug-and-play operator blocks under structure-preserving rules \cite{zhang2026legonet,zhang2026geometry}.  Our method joins the learned local subspace of GMsFEM-NO to the modular principle of LegONet, and departs from NOEM and GMsFEM-NO on three points: from NOEM in that the network reads the local coefficient alone, not the boundary state; from GMsFEM-NO in that it learns the complete port-response space of a primitive, not a spectral subset; and from both in that the composition is recursive, under a localized bound on the error the learned spaces introduce.  On our two-dimensional benchmark, NOEM reaches $10\%$ energy error against $4.4\%$ for Ms-LegONet on the same skeleton unknowns (\Cref{subsec:2d-diffusion}).

\begin{figure}[!t]
\centering
\includegraphics[width=0.99\textwidth]{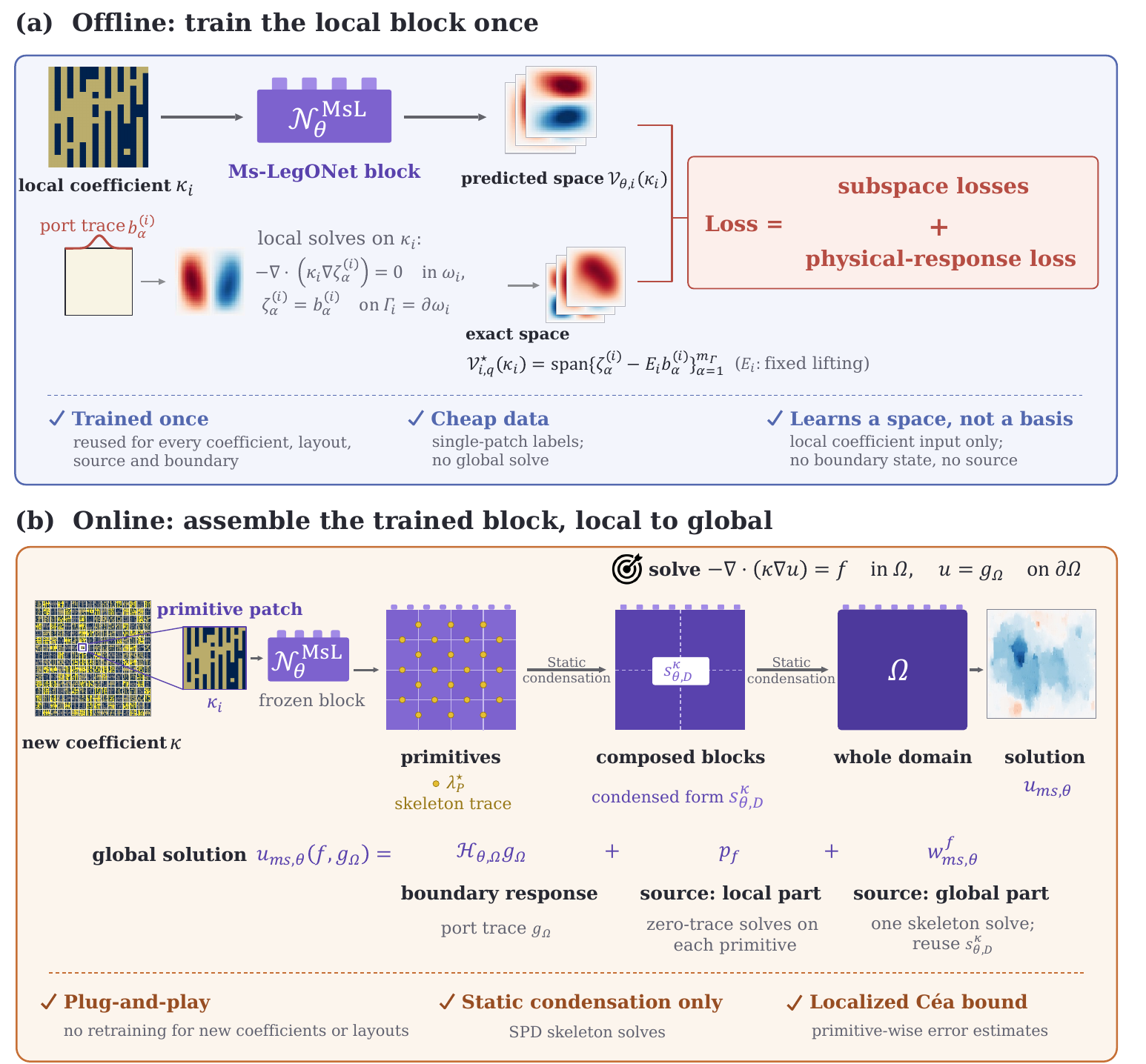}
\caption{Workflow of Ms-LegONet. (a) Offline, the block
$\mathcal N_\theta^{\mathrm{MsL}}$ is trained on single-patch zero-trace
corrections to predict the coefficient-dependent local space. (b) Online,
one evaluation per primitive gives the space and each port trace is extended
by a Galerkin solve in it; blocks are composed recursively by static
condensation up to the whole domain. Source and boundary data enter through
the three-term decomposition \eqref{eq:global-ms-decomposition}: the
boundary datum as the root port trace, local particular fields
\eqref{eq:local-particular-problem}, and one skeleton solve
\eqref{eq:global-ms-galerkin} for the source-induced correction; no further network is
introduced.}
\label{fig:workflow}
\end{figure}

We propose \emph{multiscale LegONet} (Ms-LegONet), a plug-and-play multiscale block for elliptic problems (\Cref{fig:workflow}).  On a reference primitive patch $\widehat\omega$ one network per spatial dimension learns
\begin{equation}
    \widehat\kappa
    \longmapsto
    \widehat{\mathcal V}_{\theta}(\widehat\kappa)
    =\spanop\{\psi_{\theta,1},\ldots,\psi_{\theta,q}\}
    \subset H_0^1(\widehat\omega),
    \label{eq:intro-learned-map}
\end{equation}
where the $\psi_{\theta,r}$ are physical zero-boundary fields.  Offline, the block is trained once on the zero-trace corrections of single patches, at a data cost independent of domain size and layout.  Online, one network evaluation per primitive gives the space, a fixed conforming trace basis provides the ports, and each port's harmonic field is a Galerkin solve in the learned space.  Larger blocks are built by recursive static condensation, with no network at parent scales; source and boundary data enter through local particular fields, the root port trace and one skeleton solve, so the global solution is conforming Galerkin, under a localized C\'ea bound that attributes the learned error primitive by primitive.  Ms-LegONet is thus a multiscale finite element method with a learned local space: it keeps the coarse variational structure of MsFEM and GMsFEM and the component--port composition of \cite{huynh2013scrbe}, and obtains the local basis from one network evaluation per primitive.

We make four contributions.
\begin{itemize}[leftmargin=2em]
    \item \emph{The learned object is a space, not a basis.}  The network predicts the coefficient-dependent local space \eqref{eq:intro-learned-map} of physical fields, not a stiffness, Schur or factor representation; fixing the space, not its basis, lets the downstream construction choose its port basis, sources and boundary data without retraining.
    \item \emph{Plug-and-play across coefficients.}  One network per dimension, trained once, serves every new coefficient of its class, so the block can be the inner solve of a many-query sweep.
    \item \emph{The learned block yields large speed-ups.}  Recursive static condensation assembles and parallelizes the learned blocks.
    For a 64-primitive chain in three-dimensional layered media, Ms-LegONet computes the $32\times32$ port-response matrix in 0.11 s per coefficient field, a $58\times$ speedup over $Q_1$ blocks on the same mesh, at a relative error of $10^{-4}$ against their $6\%$ (\Cref{tab:3d-eff}).
     For the global solution, a monolithic solve of the accuracy of the learned block needs about $100^3$ cells per primitive and three to four orders of magnitude more time (\Cref{tab:3d-baselines}).
    \item \emph{A localized error bound.}  The conforming Galerkin solution obeys the localized C\'ea bound of \Cref{prop:localized-bound}, verified on every two-dimensional test mosaic--problem pair with median effectivity $1.002$--$1.034$ on the mixed layouts and at most $1.093$ elsewhere, and pins the error to the primitives that cause it.
\end{itemize}

\Cref{sec:problem} states the model problem; \Cref{sec:method} develops the block, its composition, sources, boundary data and the localized bound; \Cref{sec:implementation} describes network, data and loss; \Cref{sec:numerical-results} reports the results; \Cref{sec:conclusion} concludes; the appendices contain the proofs.
\section{Problem formulation}
\label{sec:problem}

Let \(\Omega\subset\mathbb{R}^d\), \(d\in\{2,3\}\), be a bounded
connected Lipschitz domain and \(f\in H^{-1}(\Omega):=(H_0^1(\Omega))'\).
We consider the heterogeneous diffusion problem
\begin{align}
-\nabla\cdot\bigl(\kappae\nabla u\bigr)
&=f
&&\text{in }\Omega,
\label{eq:pde}
\\
u&=0
&&\text{on }\partial\Omega,
\label{eq:bc}
\end{align}
in the weak sense; nonhomogeneous Dirichlet data are treated in
\Cref{subsec:source-global}. The coefficient
\(\kappae\in L^\infty(\Omega)\) satisfies
\begin{equation}
0<\kappa_{\min}
\le
\kappae(\mathbf{x})
\le
\kappa_{\max}
<\infty
\qquad
\text{for a.e. }\mathbf{x}\in\Omega,
\label{eq:kappa-positive}
\end{equation}
with \(\kappa_{\min}\), \(\kappa_{\max}\) its essential bounds and
contrast \(\chi(\kappae):=\kappa_{\max}/\kappa_{\min}\). The
coefficient families (discontinuous, random, high-contrast
\cite{efendiev2009multiscale,efendiev2013gmsfem,chung2023multiscale})
are specified in \Cref{sec:numerical-results}.

Set \(V:=H_0^1(\Omega)\) with \(\|v\|_V:=\|\nabla v\|_{L^2(\Omega)}\)
\cite{evans2022partial,brenner2008mathematical}, and for \(w,v\in V\)
let \(a_{\kappae}(w,v):=\int_{\Omega}\kappae\,\nabla w\cdot\nabla v\,\mathrm{d}\mathbf{x}\)
and \(\ell_f(v):=\langle f,v\rangle_{H^{-1}(\Omega),H_0^1(\Omega)}\).
The weak form of \eqref{eq:pde}--\eqref{eq:bc} is: find \(u\in V\) such
that
\begin{equation}
a_{\kappae}(u,v)=\ell_f(v)
\qquad
\forall v\in V.
\label{eq:weak-problem}
\end{equation}
By \eqref{eq:kappa-positive}, \(a_{\kappae}\) is continuous with constant
\(\kappa_{\max}\) and coercive with constant \(\kappa_{\min}\) on \(V\),
so \eqref{eq:weak-problem} has a unique solution (Lax--Milgram
\cite{evans2022partial,brenner2008mathematical,ern2004theory}), and the
energy norm \(\|v\|_{a_{\kappae}}:=a_{\kappae}(v,v)^{1/2}\) satisfies
\(\sqrt{\kappa_{\min}}\,\|v\|_V\le\|v\|_{a_{\kappae}}\le\sqrt{\kappa_{\max}}\,\|v\|_V\).

The fine-scale reference uses a conforming, shape-regular tensor-product
mesh \(\mathcal{T}_h\) of \(\Omega\) (quadrilateral or hexahedral) with
\(h:=\max_{K\in\mathcal{T}_h}\operatorname{diam}(K)\). With
\(Q_1(K)\) the space of multilinear functions on \(K\), the finite element
space is
\(V_h:=\{v_h\in H_0^1(\Omega):v_h|_K\in Q_1(K)\ \forall K\in\mathcal{T}_h\}\subset V\),
and the fine-mesh solution \(u_h\in V_h\) satisfies
\begin{equation}
	a_{\kappae}(u_h,v_h)
	=
	\ell_f(v_h)
	\qquad
	\forall v_h\in V_h.
	\label{eq:fine-fe}
\end{equation}
Galerkin orthogonality and the symmetry of \(a_{\kappae}\) give
\(\|u-u_h\|_{a_{\kappae}}=\inf_{v_h\in V_h}\|u-v_h\|_{a_{\kappae}}\),
whereas C\'ea's lemma in the \(V\)-norm
\cite{brenner2008mathematical,ern2004theory} carries a constant that
grows with \(\kappa_{\max}/\kappa_{\min}\). The energy norm is therefore used
throughout: in it, every Galerkin solution below is a best approximation
with constant one, independent of the contrast.

\section{Multiscale LegONet}
\label{sec:method}

The method follows the stages of \Cref{fig:workflow}. \Cref{subsec:primitive,subsec:neural-block,subsec:composition}
treat the source-free problem, in which the learned block approximates the
\(\kappa\)-harmonic response of a primitive to its boundary trace and
recursive composition assembles these responses up to \(\Omega\); source
and boundary data enter only in \Cref{subsec:source-global}, so that
neither is an input of the network, and \Cref{subsec:localized-bound}
bounds the resulting error.

\subsection{Primitive patches and conforming boundary ports}
\label{subsec:primitive}

Let
\(\mathcal T_H=\{\omega_i\}_{i=1}^{N_H}\)
be a conforming, nonoverlapping coarse partition of \(\Omega\) into
connected Lipschitz \emph{primitive patches}, aligned with the fine
mesh, with
\(H:=\max_{i}\operatorname{diam}(\omega_i)\) and boundaries
\(\Gamma_i:=\partial\omega_i\). For adjacent primitives,
\(\Gamma_{ij}:=\partial\omega_i\cap\partial\omega_j\) with
\(|\Gamma_{ij}|_{d-1}>0\) (\(|\cdot|_{d-1}\) the
\((d-1)\)-dimensional Hausdorff measure) is an interior interface;
\(\mathcal F_H^{\mathrm{int}}\) is the set of these.

Let \(\widehat\omega\subset\mathbb R^d\) be a fixed reference primitive
with conforming reference mesh \(\widehat{\mathcal T}_h\) and
\(F_i:\widehat\omega\rightarrow\omega_i\) invertible affine maps with
\(\mathcal T_h(\omega_i)=F_i(\widehat{\mathcal T}_h)\); the affine
pull-back is detailed in \Cref{app:affine-pullback}. With \(\kappa_i:=\kappae|_{\omega_i}\),
the local energy form is
\(a_i^{\kappa_i}(w,v):=\int_{\omega_i}\kappa_i\nabla w\cdot\nabla v\,\mathrm d\mathbf x\)
for \(w,v\in H^1(\omega_i)\).
Let \(V_h(\omega_i)\) be the conforming \(Q_1\) space on
\(\mathcal T_h(\omega_i)\),
\(\gamma_i:H^1(\omega_i)\rightarrow H^{1/2}(\Gamma_i)\) the trace
operator, \(V_{i,h}^{0}:=\ker(\gamma_i|_{V_h(\omega_i)})\) the local
zero-trace space and \(\Lambda_h(\Gamma_i):=\gamma_iV_h(\omega_i)\) the
discrete trace space, with reference
counterparts \(\widehat V_h\), \(\widehat\gamma\) and
\(\widehat\Lambda_h:=\widehat\gamma\widehat V_h\). Fix a
coefficient-independent port space
\begin{equation}
\widehat{\mathcal B}_{\Gamma}
:=
\operatorname{span}
\{\widehat b_1,\ldots,\widehat b_{m_\Gamma}\}
\subseteq
\widehat\Lambda_h
\label{eq:reference-port-space}
\end{equation}
of dimension \(m_\Gamma\), with physical realization
\(\mathcal B_{\Gamma_i}:=\{\widehat g\circ F_i^{-1}:\widehat g\in\widehat{\mathcal B}_{\Gamma}\}\subseteq\Lambda_h(\Gamma_i)\)
and port traces \(b_\alpha^{(i)}:=\widehat b_\alpha\circ F_i^{-1}\).
The port spaces are assumed compatible on every
\(\Gamma_{ij}\in\mathcal F_H^{\mathrm{int}}\): the restrictions of
\(\mathcal B_{\Gamma_i}\) and \(\mathcal B_{\Gamma_j}\) to \(\Gamma_{ij}\) span
the same space once \(\Gamma_{ij}\) is parametrized consistently from both
sides, so shared interface trace coefficients enforce
\begin{equation}
\left.(\gamma_i v_i)\right|_{\Gamma_{ij}}
=
\left.(\gamma_j v_j)\right|_{\Gamma_{ij}}
\qquad
\text{a.e. on }\Gamma_{ij},
\label{eq:interface-conformity}
\end{equation}
and the piecewise field is \(H^1\)-conforming
\cite{brenner2008mathematical,ern2004theory}.

For \(g\in\mathcal B_{\Gamma_i}\), the local
\(\kappa_i\)-harmonic extension
\(\mathcal H_i^{\kappa_i}g\) minimizes \(\frac12a_i^{\kappa_i}(v,v)\)
over \(v\in H^1(\omega_i)\) with \(\gamma_iv=g\);
its conforming counterpart \(\mathcal H_{i,h}^{\kappa_i}g\) restricts
the minimization to \(V_h(\omega_i)\). Both exist, are unique and depend
linearly on \(g\) by the Lax--Milgram theorem
\cite{brenner2008mathematical,ern2004theory}. With
\(\zeta_\alpha^{(i)}
:=\mathcal H_{i,h}^{\kappa_i}b_\alpha^{(i)}\), the local
harmonic-response space is
\begin{equation}
\mathcal W_h(\omega_i;\kappa_i)
:=
\mathcal H_{i,h}^{\kappa_i}\mathcal B_{\Gamma_i}
=
\operatorname{span}
\left\{
\zeta_\alpha^{(i)}
\right\}_{\alpha=1}^{m_\Gamma}
\subset V_h(\omega_i).
\label{eq:local-harmonic-response-space}
\end{equation}
Assembling \(\mathcal W_h(\omega_i;\kappa_i)\) costs \(m_\Gamma\) local
solves, which MsFEM and GMsFEM
\cite{hou1997msfem,efendiev2011spectral,efendiev2013gmsfem} repeat for
every new coefficient \(\kappa_i\); Ms-LegONet replaces them by a learned
coefficient-dependent space.

\subsection{Local Ms-LegONet block}
\label{subsec:neural-block}

Fix a coefficient-independent linear lifting
\(E_i:\mathcal B_{\Gamma_i}\rightarrow V_h(\omega_i)\) with
\(\gamma_iE_i g=g\) for every \(g\in\mathcal B_{\Gamma_i}\) (one
exists because
\(\mathcal B_{\Gamma_i}\subseteq\gamma_iV_h(\omega_i)\)); our
two-dimensional experiments use a Coons lift \cite{coons1967surfaces,gordon1973transfinite}. The
zero-trace correction operator is
\begin{equation}
\mathcal C_{i,h}^{\kappa_i}
:=
\mathcal H_{i,h}^{\kappa_i}-E_i,
\qquad
\mathcal H_{i,h}^{\kappa_i}g
=
E_i g+\mathcal C_{i,h}^{\kappa_i}g,
\qquad
\mathcal C_{i,h}^{\kappa_i}g\in V_{i,h}^{0}.
\label{eq:local-correction-operator}
\end{equation}
Let
\begin{equation}
\mathcal R_{i,h}(\kappa_i):=\operatorname{Range}(\mathcal C_{i,h}^{\kappa_i}),
\qquad
q:=\dim\mathcal R_{i,h}(\kappa_i),
\qquad
\mathcal V_{i,q}^{\star}(\kappa_i):=\mathcal R_{i,h}(\kappa_i).
\label{eq:target-local-physical-space}
\end{equation}
We assume that \(q\) does not depend on \(\kappa_i\) within the admissible
local coefficient class. The local block is trained to approximate
\(\mathcal V_{i,q}^{\star}(\kappa_i)\).

\begin{remark}
Since \(q=m_\Gamma-\dim\ker\mathcal C_{i,h}^{\kappa_i}\) and the kernel
consists of the port data whose lift \(E_ig\) is already
\(\kappa_i\)-harmonic, the assumption holds when the dimension of this
kernel is the same for every coefficient of the class. For the Coons lift
of \Cref{subsec:2d-diffusion} the kernel is generically spanned by the constant trace,
so \(q=m_\Gamma-1\); for the lift of \Cref{subsec:3d-layered} it is
trivial, so \(q=m_\Gamma\). If \(q\) varied, it would be set to its maximum
over the class, and the estimates below would hold verbatim.
\end{remark}

Pullback through \(F_i\) identifies every local zero-trace space with
\(\widehat V_h^0:=\ker(\widehat\gamma|_{\widehat V_h})\).
Let \(\widehat\kappa_i:=\kappa_i\circ F_i\) and
\(\mathcal K_{\mathrm{loc}}\) the class of reference-patch
coefficients seen in training. The shared local
block is the Grassmann-valued map
\begin{equation}
\mathcal N_\theta^{\mathrm{MsL}}:
\mathcal K_{\mathrm{loc}}
\longrightarrow
\operatorname{Gr}(q,\widehat V_h^0),
\qquad
\widehat\kappa_i
\longmapsto
\widehat{\mathcal V}_{\theta}(\widehat\kappa_i)
\approx
\widehat{\mathcal V}_{i,q}^{\star}(\widehat\kappa_i),
\label{eq:mslegonet-space-map}
\end{equation}
with \(\operatorname{Gr}(q,\widehat V_h^0)\) the Grassmann
manifold of \(q\)-dimensional subspaces of \(\widehat V_h^0\) and
\(\widehat{\mathcal V}_{i,q}^{\star}\) the pullback of
\(\mathcal V_{i,q}^{\star}(\kappa_i)\); the predicted space on
\(\omega_i\) is the push-forward
\(\mathcal V_{\theta,i}(\kappa_i)
:=
\left\{
\widehat v\circ F_i^{-1}:
\widehat v\in
\widehat{\mathcal V}_{\theta}(\widehat\kappa_i)
\right\}
\approx
\mathcal V_{i,q}^{\star}(\kappa_i)\).
One network per spatial dimension is trained once
(\Cref{sec:implementation}) and serves every primitive, layout, source,
boundary datum and composition depth.

\begin{remark}
The experiments of \Cref{sec:numerical-results} use translations, \(\mathbf B_i=\mathbf I\),
so the pulled-back coefficient of \eqref{eq:affine-pulled-coefficient} is
\(\widehat{\mathbf K}_i=\widehat\kappa_i\mathbf I\) and the scalar
\(\widehat\kappa_i\) is the complete input; for a general affine \(F_i\)
the block would read the tensor \(\widehat{\mathbf K}_i\) instead, with
\(d(d+1)/2\) input channels.
\end{remark}

Let \(\{\widehat\varphi_{k}\}_{k=1}^{n_0}\) be a basis of
\(\widehat V_h^0\), \(n_0:=\dim\widehat V_h^0\), and
\(\varphi_{k}^{(i)}:=\widehat\varphi_{k}\circ F_i^{-1}\).
The block returns a full-column-rank matrix
\(\widetilde{\mathbf Q}_{\theta,i}\in\mathbb R^{n_0\times q}\)
representing
\begin{equation}
\mathcal V_{\theta,i}(\kappa_i)
=
\operatorname{span}
\left\{
\widetilde\psi_{\theta,r}^{(i)}
:=
\sum_{k=1}^{n_0}
\bigl(\widetilde{\mathbf Q}_{\theta,i}\bigr)_{kr}
\varphi_{k}^{(i)}
\right\}_{r=1}^{q}.
\label{eq:neural-physical-space}
\end{equation}
With the interior stiffness matrix
\(\mathbf A_i(\kappa_i)\in\mathbb R^{n_0\times n_0}\), define
\begin{equation}
\begin{aligned}
\bigl[\mathbf A_i(\kappa_i)\bigr]_{kl}
&:=
a_i^{\kappa_i}
\bigl(\varphi_{l}^{(i)},\varphi_{k}^{(i)}\bigr),\\
\mathbf G_{\theta,i}(\kappa_i)
&:=
\widetilde{\mathbf Q}_{\theta,i}^{T}
\mathbf A_i(\kappa_i)
\widetilde{\mathbf Q}_{\theta,i},\\
\mathbf Q_{\theta,i}
&:=
\widetilde{\mathbf Q}_{\theta,i}
\mathbf G_{\theta,i}^{-1/2},
\end{aligned}
\label{eq:energy-orthonormalized-neural-basis}
\end{equation}
where \(\mathbf G_{\theta,i}^{-1/2}\) is the symmetric positive-definite
inverse square root. For any
\(\mathbf A_i(\kappa_i)\)-orthonormal basis
\(\mathbf Q_i^\star\) of
\(\mathcal V_{i,q}^{\star}(\kappa_i)\), the energy
chordal discrepancy between the two spaces is
\(d_{\mathrm{ch}}^{2}=q-\|{\mathbf Q_i^\star}^{T}\mathbf A_i(\kappa_i)\mathbf Q_{\theta,i}\|_F^2\),
the sum of the squared sines of their energy principal angles
(\Cref{prop:app-energy-subspace-invariance}); the training loss uses
the same angles (\Cref{sec:implementation}).

For \(g\in\mathcal B_{\Gamma_i}\), let
\(u_{\theta,i}(g)=E_i g+w_{\theta,i}(g)\) with
\(w_{\theta,i}(g)\in\mathcal V_{\theta,i}(\kappa_i)\) the unique
solution of
\begin{equation}
a_i^{\kappa_i}
\bigl(
E_i g+w_{\theta,i}(g),v_\theta
\bigr)
=
0
\qquad
\forall v_\theta\in\mathcal V_{\theta,i}(\kappa_i).
\label{eq:mslegonet-local-galerkin}
\end{equation}
This is the energy projection of
\(\mathcal C_{i,h}^{\kappa_i}g\) onto the predicted space, so
\begin{equation}
\left\|
\mathcal H_{i,h}^{\kappa_i}g-u_{\theta,i}(g)
\right\|_{a_i^{\kappa_i}}
=
\inf_{v_\theta\in\mathcal V_{\theta,i}(\kappa_i)}
\left\|
\mathcal C_{i,h}^{\kappa_i}g-v_\theta
\right\|_{a_i^{\kappa_i}},
\label{eq:local-response-best-approximation}
\end{equation}
with
\(\|v\|_{a_i^{\kappa_i}}
:=a_i^{\kappa_i}(v,v)^{1/2}\), a seminorm on \(H^1(\omega_i)\) and a norm
on \(V_{i,h}^0\); recovery is exact
for every admissible port trace when
\(\mathcal V_{\theta,i}(\kappa_i)
=\mathcal V_{i,q}^{\star}(\kappa_i)\)
(\Cref{prop:app-local-response-projection}).
The network output does not depend
on \(g\): one neural evaluation per primitive fixes the space, and each
admissible trace then costs only the \(q\times q\) system
\eqref{eq:mslegonet-local-galerkin}. The learned object is a
subspace, not a solution map.

\subsection{Deterministic recursive composition}
\label{subsec:composition}

Recursive composition replaces neighboring \emph{child blocks} by one
\emph{parent block} with the same boundary-response interface; level-zero
children are the primitives, and each parent is a child at the next level.
For a primitive or composed block \(D\) with \(\Gamma_D:=\partial D\) and
port space \(\mathcal B_{\Gamma_D}\), let
\(\mathcal H_{\theta,D}:\mathcal B_{\Gamma_D}\rightarrow H^1(D)\) be its
reduced harmonic-response map: \(\mathcal H_{\theta,\omega_i}g:=u_{\theta,i}(g)\)
from \eqref{eq:mslegonet-local-galerkin} for a primitive, and
\eqref{eq:parent-response-map} below for a composed block. The energy form
on \(D\) is
\(a_D^{\kappae}(w,v):=\sum_{\omega_i\subseteq D}a_i^{\kappa_i}(w|_{\omega_i},v|_{\omega_i})\)
and the condensed form
\(s_{\theta,D}^{\kappae}(g,g'):=a_D^{\kappae}(\mathcal H_{\theta,D}g,\mathcal H_{\theta,D}g')\).

Let \(P\) be a connected parent with a finite set
\(\mathcal D(P)\) of nonoverlapping children,
\(\overline{\omega_P}=\bigcup_{D\in\mathcal D(P)}\overline D\),
external boundary \(\Gamma_P:=\partial\omega_P\) and interior skeleton
\(\Sigma_P:=(\bigcup_{D\in\mathcal D(P)}\partial D)\setminus\Gamma_P\).
Let \(\mathcal B_{\Gamma_P}\) be the external port space and
\(\Lambda_P^{\mathrm{int}}\) the finite-dimensional conforming trace
space on \(\Sigma_P\). For each child \(D\in\mathcal D(P)\),
restriction and orientation maps
\(\mathcal M_{P\to D}^{E}:
\mathcal B_{\Gamma_P}\rightarrow\mathcal B_{\Gamma_D}\)
and
\(\mathcal M_{P\to D}^{I}:
\Lambda_P^{\mathrm{int}}\rightarrow\mathcal B_{\Gamma_D}\)
map a parent trace
\(g_P\in\mathcal B_{\Gamma_P}\) and a skeleton trace
\(\lambda\in\Lambda_P^{\mathrm{int}}\) to the trace on \(D\),
\begin{equation}
g_D(g_P,\lambda)
=
\mathcal M_{P\to D}^{E}g_P
+
\mathcal M_{P\to D}^{I}\lambda.
\label{eq:child-induced-trace}
\end{equation}
The skeleton trace minimizes the total child energy:
\begin{equation}
\lambda_P^\star(g_P)
=
\operatorname*{arg\,min}_{\lambda\in\Lambda_P^{\mathrm{int}}}
J_P(\lambda;g_P),
\quad
J_P(\lambda;g_P)
:=
\frac12
\sum_{D\in\mathcal D(P)}
a_D^{\kappae}
\bigl(
\mathcal H_{\theta,D}g_D,
\mathcal H_{\theta,D}g_D
\bigr),
\label{eq:parent-energy-minimization}
\end{equation}
with \(g_D=g_D(g_P,\lambda)\) in \(J_P\); under the conformity and
injectivity assumptions of \Cref{prop:app-parent-composition} the
minimizer is unique, and the parent response map is
\begin{equation}
\left.
\mathcal H_{\theta,P}g_P
\right|_D
:=
\mathcal H_{\theta,D}
g_D\bigl(g_P,\lambda_P^\star(g_P)\bigr),
\qquad
D\in\mathcal D(P).
\label{eq:parent-response-map}
\end{equation}
It belongs to \(H^1(\omega_P)\), has trace \(g_P\) on
\(\Gamma_P\) and satisfies the reduced interface-equilibrium condition
\begin{equation}
\sum_{D\in\mathcal D(P)}
a_D^{\kappae}
\left(
\mathcal H_{\theta,D}
g_D\bigl(g_P,\lambda_P^\star\bigr),
\mathcal H_{\theta,D}
\mathcal M_{P\to D}^{I}\mu
\right)
=
0
\qquad
\forall\mu\in\Lambda_P^{\mathrm{int}}.
\label{eq:interface-equilibrium}
\end{equation}
With the condensed form
\(s_{\theta,P}^{\kappae}(g,g')
:=a_P^{\kappae}
\bigl(\mathcal H_{\theta,P}g,\mathcal H_{\theta,P}g'\bigr)\),
the parent can be composed again without another neural model. The
minimization \eqref{eq:parent-energy-minimization} is a symmetric
positive-definite skeleton solve (static condensation
\cite{quarteroni1999domain,toselli2004domain}; matrix form in
\Cref{app:recursive-composition-algebra}).

Let \(\{\mathcal T_H^{(j)}\}_{j=0}^{L}\) be a nested block
hierarchy with \(\mathcal T_H^{(0)}=\mathcal T_H\) and
\(\mathcal T_H^{(L)}=\{\Omega\}\), each level built from the previous
by
\eqref{eq:parent-energy-minimization}--%
\eqref{eq:parent-response-map}; neural inference occurs only at level
zero. We take the parent port space \(\mathcal B_{\Gamma_P}\) to be the
restriction of the children's port traces to \(\Gamma_P\) and
\(\Lambda_P^{\mathrm{int}}\) to be the space of all conforming interior
trace tuples, so that no trace is truncated at a parent level; both
conditions hold for the nodal traces of \Cref{subsec:2d-diffusion} and the chains of
\Cref{subsec:3d-layered}. Minimizing the total energy level by level is then the same as
minimizing it over the whole level-zero skeleton at once: nested and
single-shot elimination of a symmetric positive-definite system give the
same Schur complement on the retained unknowns
(\Cref{app:recursive-composition-algebra}). The hierarchy changes the
computation, not the approximation.

\subsection{Global problem with forcing and boundary data}
\label{subsec:source-global}

For nonhomogeneous boundary data
\(g_\Omega\in H^{1/2}(\partial\Omega)\) consider
\(-\nabla\cdot(\kappae\nabla u)=f\) in \(\Omega\) with \(u=g_\Omega\) on
\(\partial\Omega\), of which \eqref{eq:pde}--\eqref{eq:bc} is the case
\(g_\Omega=0\). Let
\(\overline V_h:=\{v_h\in H^1(\Omega):v_h|_K\in Q_1(K)\ \forall K\in\mathcal T_h\}\)
be the unconstrained \(Q_1\) space, \(\gamma\) the trace operator on
\(\partial\Omega\), so that \(V_h=\ker(\gamma|_{\overline V_h})\), and
\(\overline V_h^{g_\Omega}:=\{v_h\in\overline V_h:\gamma v_h=g_\Omega\}\)
for a discrete datum \(g_\Omega\in\gamma\overline V_h\). The fine-grid
reference solution \(u_h(f,g_\Omega)\in\overline V_h^{g_\Omega}\)
satisfies
\begin{equation}
a_{\kappae}\bigl(u_h(f,g_\Omega),v_h\bigr)
=
\ell_f(v_h)
\qquad
\forall v_h\in V_h,
\label{eq:global-source-discrete-weak-form}
\end{equation}
and reduces to \(u_h\) of \eqref{eq:fine-fe} for \(g_\Omega=0\).

For \(g_\Omega\in\mathcal B_{\Gamma_\Omega}
\subseteq\gamma\overline V_h\), with
\(\Gamma_\Omega:=\partial\Omega\) and \(\mathcal B_{\Gamma_\Omega}\)
the port space of the root block, we seek
\begin{equation}
u_{\mathrm{ms},\theta}(f,g_\Omega)
=
u_{\mathrm{ms},\theta}^{g_\Omega}
+
p_f
+
w_{\mathrm{ms},\theta}^{f},
\label{eq:global-ms-decomposition}
\end{equation}
with boundary-driven harmonic response
\(u_{\mathrm{ms},\theta}^{g_\Omega}\), local particular
contribution \(p_f\) and source-induced harmonic correction
\(w_{\mathrm{ms},\theta}^{f}\); no additional network is introduced
for \(f\) or \(g_\Omega\).

\paragraph{Boundary-driven harmonic response}
The first term is the root response
\begin{equation}
u_{\mathrm{ms},\theta}^{g_\Omega}
:=
\mathcal H_{\theta,\Omega}g_\Omega,
\label{eq:boundary-harmonic-contribution}
\end{equation}
with \(\mathcal H_{\theta,\Omega}\) the recursive composition
\eqref{eq:parent-energy-minimization}--\eqref{eq:parent-response-map}
up to the root block.

\paragraph{Local particular contribution}
For each primitive \(\omega_i\) let
\(\mathcal E_i^0:V_{i,h}^0\rightarrow V_h\) be extension by zero and
\(\ell_{f,i}(v_i):=\ell_f(\mathcal E_i^0v_i)\) the restricted load.
The local particular field \(p_i(f)\in V_{i,h}^0\) and its
assembly \(p_f\in V_h\) are defined by
\begin{equation}
\begin{aligned}
a_i^{\kappa_i}\bigl(p_i(f),v_i\bigr)
&=
\ell_{f,i}(v_i)
&&
\forall v_i\in V_{i,h}^0,
\\
p_f
&:=
\sum_{i=1}^{N_H}\mathcal E_i^0p_i(f).
\end{aligned}
\label{eq:local-particular-problem}
\end{equation}
These local problems contain no learned quantity. By
\Cref{prop:app-particular-harmonic-decomposition},
\(u_h(f,g_\Omega)-p_f\) is discretely \(\kappa_i\)-harmonic on every
primitive with external trace \(g_\Omega\), hence approximable in the
composed harmonic-response spaces.

\paragraph{Global harmonic equilibration}
Let \(\Lambda_{H,0}^{\mathrm{sk}}\) be the conforming zero-boundary
skeleton space: the tuples
\((\lambda_i)_{i=1}^{N_H}\in\prod_{i=1}^{N_H}\mathcal B_{\Gamma_i}\) that agree
on every interior interface in the sense of \eqref{eq:interface-conformity}
and vanish on \(\Gamma_i\cap\partial\Omega\). The global
zero-boundary multiscale space is
\begin{equation}
V_{\mathrm{ms},\theta}^{0}
:=
\left\{
v_\lambda\in V_h:
\left.v_\lambda\right|_{\omega_i}
=
\mathcal H_{\theta,\omega_i}\lambda_i
\quad\forall i,
\quad
(\lambda_i)_{i=1}^{N_H}
\in\Lambda_{H,0}^{\mathrm{sk}}
\right\},
\label{eq:global-zero-boundary-ms-space}
\end{equation}
and recursive elimination of the skeleton variables is static
condensation in \(V_{\mathrm{ms},\theta}^{0}\)
\cite{quarteroni1999domain,toselli2004domain}.
Because \(u_{\mathrm{ms},\theta}^{g_\Omega}\) is the recursively
energy-minimizing root response,
\begin{equation}
a_{\kappae}
\bigl(
u_{\mathrm{ms},\theta}^{g_\Omega},
v_{\mathrm{ms},\theta}
\bigr)
=
0
\qquad
\forall
v_{\mathrm{ms},\theta}\in V_{\mathrm{ms},\theta}^{0}.
\label{eq:boundary-response-orthogonality}
\end{equation}
The source-induced correction is the unique
\(w_{\mathrm{ms},\theta}^{f}\in V_{\mathrm{ms},\theta}^{0}\) with
\begin{equation}
a_{\kappae}
\bigl(
w_{\mathrm{ms},\theta}^{f},
v_{\mathrm{ms},\theta}
\bigr)
=
\ell_f(v_{\mathrm{ms},\theta})
-
a_{\kappae}
\bigl(
p_f,v_{\mathrm{ms},\theta}
\bigr)
\qquad
\forall
v_{\mathrm{ms},\theta}
\in V_{\mathrm{ms},\theta}^{0};
\label{eq:global-ms-galerkin}
\end{equation}
by
\eqref{eq:global-ms-decomposition},
\eqref{eq:boundary-response-orthogonality} and
\eqref{eq:global-ms-galerkin},
\(a_{\kappae}(u_{\mathrm{ms},\theta}(f,g_\Omega),
v_{\mathrm{ms},\theta})
=\ell_f(v_{\mathrm{ms},\theta})\)
for every \(v_{\mathrm{ms},\theta}\in V_{\mathrm{ms},\theta}^{0}\),
so the Ms-LegONet solution is the Galerkin solution in the conforming
affine space
\(u_{\mathrm{ms},\theta}^{g_\Omega}+p_f+V_{\mathrm{ms},\theta}^{0}\).
For \(g_\Omega=0\),
\(u_{\mathrm{ms},\theta}^{g_\Omega}=0\)
and
\(u_{\mathrm{ms},\theta}(f,0)=p_f+w_{\mathrm{ms},\theta}^{f}\).
\Cref{alg:mslegonet} summarizes the method: the network is evaluated
once per primitive, and every later step is deterministic linear algebra,
including the local solves for \(p_i(f)\) when \(f\neq0\) (line 5).

\begin{algorithm}[!htbp]
\caption{Ms-LegONet}
\label{alg:mslegonet}
\begin{algorithmic}[1]
\Require coefficient \(\kappae\), source \(f\), boundary datum \(g_\Omega\), partition \(\mathcal T_H\) with hierarchy \(\{\mathcal T_H^{(j)}\}_{j=0}^{L}\), port basis \(\{\widehat b_\alpha\}\) \eqref{eq:reference-port-space}, lifting \(E_i\)
\Ensure \(u_{\mathrm{ms},\theta}(f,g_\Omega)\) of \eqref{eq:global-ms-decomposition}
\Statex \emph{Offline, once per spatial dimension:}
\State train \(\theta\) on single-patch harmonic responses by minimizing \eqref{eq:implementation-total-loss}; freeze \(\theta\)
\Statex \emph{Online, for each new \((\kappae,f,g_\Omega)\):}
\For{each primitive \(\omega_i\in\mathcal T_H\), in parallel}
  \State \(\widetilde{\mathbf Q}_{\theta,i}\gets\mathcal N_\theta^{\mathrm{MsL}}(\widehat\kappa_i)\) \eqref{eq:mslegonet-space-map}; orthonormalize by \eqref{eq:energy-orthonormalized-neural-basis}
  \State port responses \(u_{\theta,i}(b_\alpha^{(i)})\), \(\alpha=1,\ldots,m_\Gamma\), by \eqref{eq:mslegonet-local-galerkin}, explicit in the orthonormal basis of line 3; condensed form \(s_{\theta,\omega_i}^{\kappae}\)
  \State particular field \(p_i(f)\) by \eqref{eq:local-particular-problem}
\EndFor
\For{levels \(j=1,\ldots,L\)}
  \State for each parent \(P\in\mathcal T_H^{(j)}\): skeleton solve \eqref{eq:parent-energy-minimization}, response map \eqref{eq:parent-response-map}, condensed form \(s_{\theta,P}^{\kappae}\)
\EndFor
\State \(u_{\mathrm{ms},\theta}^{g_\Omega}\gets\mathcal H_{\theta,\Omega}g_\Omega\) \eqref{eq:boundary-harmonic-contribution}; \(p_f\gets\sum_i\mathcal E_i^0p_i(f)\)
\State \(w_{\mathrm{ms},\theta}^{f}\gets\) solution of \eqref{eq:global-ms-galerkin} in \(V_{\mathrm{ms},\theta}^{0}\)
\State \Return \(u_{\mathrm{ms},\theta}^{g_\Omega}+p_f+w_{\mathrm{ms},\theta}^{f}\)
\end{algorithmic}
\end{algorithm}

\subsection{A localized bound for the learned error}
\label{subsec:localized-bound}

Let \(u_{\mathrm{ms}}^{\star}\)
be the \emph{exact-block} solution: the construction of
\Cref{subsec:composition,subsec:source-global} with the exact
spaces \(\mathcal V_{i,q}^{\star}(\kappa_i)\) in place of
\(\mathcal V_{\theta,i}(\kappa_i)\), for the same port space, source
and boundary datum. Its zero-trace correction on \(\omega_i\) is
\(w_i^{\star}:=\mathcal C_{i,h}^{\kappa_i}g_i^{\star}\), with
\(g_i^{\star}\) the trace of \(u_{\mathrm{ms}}^{\star}\) on
\(\Gamma_i\). Let \(\Pi_{\theta,i}\) be the
\(a_i^{\kappa_i}\)-orthogonal projection of \(V_{i,h}^0\) onto
\(\mathcal V_{\theta,i}(\kappa_i)\).

\begin{proposition}[Localized bound]
\label{prop:localized-bound}
Let \(u_h\) be the fine-grid solution \eqref{eq:global-source-discrete-weak-form} and \(u_{\mathrm{ms},\theta}\) the Ms-LegONet solution \eqref{eq:global-ms-decomposition}, and let
\begin{equation}
\begin{aligned}
\eta_{\mathrm{proj}}
&:=
\Bigl(\sum_{i=1}^{N_H}\|(I-\Pi_{\theta,i})w_i^{\star}\|_{a_i^{\kappa_i}}^2\Bigr)^{1/2},
\\
\eta_{\mathrm{ang}}
&:=
\Bigl(\sum_{i=1}^{N_H}\sin^2\vartheta_{\max,i}\,\|w_i^{\star}\|_{a_i^{\kappa_i}}^2\Bigr)^{1/2},
\end{aligned}
\label{eq:localized-bound-terms}
\end{equation}
where \(\vartheta_{\max,i}\) is the largest energy principal angle between \(\mathcal V_{i,q}^{\star}(\kappa_i)\) and \(\mathcal V_{\theta,i}(\kappa_i)\). Then
\begin{equation}
\|u_h-u_{\mathrm{ms},\theta}\|_{a_{\kappae}}
\le
\|u_h-u_{\mathrm{ms}}^{\star}\|_{a_{\kappae}}+\eta_{\mathrm{proj}}
\le
\|u_h-u_{\mathrm{ms}}^{\star}\|_{a_{\kappae}}+\eta_{\mathrm{ang}}.
\label{eq:localized-bound}
\end{equation}
The term \(\eta_{\mathrm{proj}}\) is the best-approximation error of the exact local corrections \(w_i^{\star}\) in the learned spaces; \(\eta_{\mathrm{ang}}\) separates the accuracy of the learned space, \(\vartheta_{\max,i}\), from the size of the correction, \(\|w_i^{\star}\|_{a_i^{\kappa_i}}\). If the port space is the complete trace space \(\Lambda_h(\Gamma_i)\) on every primitive, \(u_{\mathrm{ms}}^{\star}=u_h\) and the first term vanishes.
\end{proposition}

The proof is in \Cref{app:localized-bound-proof}.
Both bounds are independent of the composition depth. They require
\(u_{\mathrm{ms}}^{\star}\), and hence exact local solves on every
primitive, so they serve to analyze the learned error, not as a posteriori
estimators (\Cref{sec:numerical-results}). The constructions and
estimates of this section use only the finite dimension of the port and
learned spaces, the \(H^1\)-conformity of the piecewise fields and Galerkin
orthogonality in the energy inner product; they hold verbatim with
\(V_h(\omega_i)\) replaced by \(H^1(\omega_i)\) and \(u_h\) by the continuum
solution \(u\), which is the setting of \Cref{subsec:3d-layered}.
\section{Implementation details}
\label{sec:implementation}

The network realizes the coefficient-to-space map
\eqref{eq:mslegonet-space-map}; every other step of
\Cref{sec:method} is deterministic.

\paragraph{Network input and output}
The input \(\mathbf X_i\) is the cellwise coefficient of the reference
primitive in normalized form,
\(z_{i,\boldsymbol\nu}:=(\log\kappa_{i,\boldsymbol\nu}-\overline{\log\kappa_i})/s_\kappa\),
with \(\kappa_{i,\boldsymbol\nu}\) the value of \(\widehat\kappa_i\) on the
cell \(\boldsymbol\nu\) of \(\widehat{\mathcal T}_h\),
\(\overline{\log\kappa_i}\) its cell mean and \(s_\kappa>0\) a fixed scale
of the coefficient class. It is concatenated with the $d$ channels that hold the reference coordinates of the cell centres; neither a boundary trace nor a source term is supplied. With
\(n_0=\dim\widehat V_h^{0}\) and \(q=\dim\mathcal R_{i,h}(\kappa_i)\), the
map \eqref{eq:mslegonet-space-map} is implemented as
\begin{equation}
\mathbf X_i
\longmapsto
\widetilde{\mathbf Q}_{\theta,i}(\kappa_i)
\in\mathbb R^{n_0\times q}
\longmapsto
\widehat{\mathcal V}_{\theta}(\widehat\kappa_i)
=
\operatorname{Range}
\bigl(\mathbf Q_{\theta,i}(\kappa_i)\bigr)
\in
\operatorname{Gr}\bigl(q,\widehat V_h^{0}\bigr).
\label{eq:implementation-network-map}
\end{equation}
The first map is the forward pass, a convolutional encoder and a decoder
that return \(q\) zero-trace fields: nodal coefficients in
\(\widehat V_h^0\), or any parametrization of zero-trace fields whose
energy Gram matrix \(\mathbf G_{\theta,i}\) can be evaluated
(\Cref{subsec:3d-layered} uses a low-rank separated representation). The
second map is the orthonormalization
\eqref{eq:energy-orthonormalized-neural-basis}; Gram matrix,
orthonormalization and quadratic forms are evaluated in higher precision
than the forward pass and remain differentiable. The target is
\(\mathcal R_{i,h}(\kappa_i)\) of \eqref{eq:target-local-physical-space},
with \(\mathbf Q_i^\star(\kappa_i)\in\mathbb R^{n_0\times q}\) an
\(\mathbf A_i(\kappa_i)\)-orthonormal basis of it.

\paragraph{Training data}
A label is the set of \(m_\Gamma\) zero-trace corrections
\(\mathcal C_{i,h}^{\kappa_i}b_\alpha^{(i)}\) of one primitive patch, whose
span is the target \(\mathcal R_{i,h}(\kappa_i)\); each is one local solve
on that patch, and no global solution enters the data. Patches are split
into disjoint training, validation and test sets, and the mosaics and
chains of \Cref{sec:numerical-results} are assembled from the
corresponding split.

\paragraph{Loss}
For each realization \(\kappa_i\), the singular values
\(\{\sigma_{i,r}\}_{r=1}^{q}\) of
\(\mathbf C_{\theta,i}:={\mathbf Q_i^\star}^{T}\mathbf A_i(\kappa_i)\mathbf Q_{\theta,i}\)
are the cosines of the energy principal angles between target and
predicted space (\Cref{prop:app-energy-subspace-invariance},
\cite{bjorck1973numerical,absil2008optimization}). With the defects
\(\delta_{i,r}:=1-\sigma_{i,r}^{2}\), ordered as
\(\delta_{i,[1]}\ge\cdots\ge \delta_{i,[q]}\), the subspace losses are
\begin{equation}
\mathcal L_{\mathrm{space},i}
:=
\frac{1}{q}
\sum_{r=1}^{q}\delta_{i,r},
\qquad
\mathcal L_{\mathrm{tail},i}
:=
\frac{1}{q_{\mathrm{tail}}}
\sum_{r=1}^{q_{\mathrm{tail}}}\delta_{i,[r]},
\label{eq:implementation-subspace-loss}
\end{equation}
with \(q_{\mathrm{tail}}\in\{1,\ldots,q\}\) prescribed.
A physical-response loss tests the predicted space on \(M\) fixed
probe traces \(g_{i,m}\in\mathcal B_{\Gamma_i}\), drawn once as random
combinations of the port basis. With
	\(e_{i,m}:=\|\mathcal H_{i,h}^{\kappa_i}g_{i,m}-u_{\theta,i}(g_{i,m})\|_{a_i^{\kappa_i}}
	/\|\mathcal H_{i,h}^{\kappa_i}g_{i,m}\|_{a_i^{\kappa_i}}\) the relative
	energy error of the local response \eqref{eq:mslegonet-local-galerkin},
	\(e_{i,[1]}\ge\cdots\ge e_{i,[M]}\) its order statistics and
	\(e^{L^2}_{i,m}\) the same ratio in the \(\kappa_i\)-weighted \(L^2\)
	norm,
	\begin{equation}
		\mathcal L_{\mathrm{phys},i}
		:=
		\frac1M\sum_{m=1}^{M}e_{i,m}
		+\frac{1}{\lceil\rho M\rceil}\sum_{m=1}^{\lceil\rho M\rceil}e_{i,[m]}
		+\frac1M\sum_{m=1}^{M}e^{L^2}_{i,m},
		\qquad \rho\in(0,1],
		\label{eq:implementation-phys-loss}
	\end{equation}
	the middle term being the conditional value at risk of the probe
	errors \cite{rockafellar2000optimization}. Like the subspace losses,
it depends on \(\mathcal V_{\theta,i}(\kappa_i)\) only, not on its
basis. With
\(\mathcal I_t\) the minibatch of update \(t\), the objective is
\begin{equation}
\mathcal J_t(\theta)
:=
\frac{1}{|\mathcal I_t|}
\sum_{i\in\mathcal I_t}
\left[
\beta_{\mathrm{space}}(t)
\mathcal L_{\mathrm{space},i}
+
\beta_{\mathrm{tail}}(t)
\mathcal L_{\mathrm{tail},i}
+
\beta_{\mathrm{phys}}(t)
\mathcal L_{\mathrm{phys},i}
\right],
\label{eq:implementation-total-loss}
\end{equation}
with nonnegative weights. Training starts from the subspace losses
and raises the physical-response weight continuously.
Optimization uses Adam \cite{kingma2014adam}; checkpoints are selected
on validation data only and frozen before testing. Networks, loss
parameters and optimization settings are listed in
\Cref{app:training-settings}.

\section{Numerical results}
\label{sec:numerical-results}

Throughout, the \emph{exact-block} solution \(u_{\mathrm{ms}}^{\star}\)
uses the exact response space \(\mathcal V_{i,q}^{\star}(\kappa_i)\) instead
of \(\mathcal V_{\theta,i}(\kappa_i)\) and isolates the effect of the
learned spaces. Every error of an approximation \(\tilde u\) to a reference
\(u_{\mathrm{ref}}\) is relative: the energy error is
\(\|u_{\mathrm{ref}}-\tilde u\|_{a_{\kappae}}/\|u_{\mathrm{ref}}\|_{a_{\kappae}}\)
and the \(L^2\) error is
\(\|u_{\mathrm{ref}}-\tilde u\|_{L^2(\Omega)}/\|u_{\mathrm{ref}}\|_{L^2(\Omega)}\).
P50 and P95 denote the median and the 95th percentile over a population.
Both \(u_{\mathrm{ms}}^{\star}\) and \(u_{\mathrm{ms},\theta}\) are Galerkin projections onto conforming spaces, and the error splits as
\begin{equation}
\|u_{\mathrm{ref}}-u_{\mathrm{ms},\theta}\|_{a_{\kappae}}
\;\le\;
\underbrace{\|u_{\mathrm{ref}}-u_{\mathrm{ms}}^{\star}\|_{a_{\kappae}}}_{\text{port truncation}}
+
\underbrace{\|u_{\mathrm{ms}}^{\star}-u_{\mathrm{ms},\theta}\|_{a_{\kappae}}}_{\text{learned space}}.
\label{eq:num-error-split}
\end{equation}
All timings are
wall-clock medians of five runs after a warm-up on a single node
(AMD EPYC 7763, 16 threads; NVIDIA A100-SXM4-40GB). Numbers quoted
without a table reference are medians over the same test sets and runs.

\subsection{Two-dimensional heterogeneous diffusion}
\label{subsec:2d-diffusion}

\paragraph{Primitive, ports, learned space and coefficient population}
The reference primitive is \(\widehat\omega=(0,1)^2\) with \(16\times16\)
bilinear \(Q_1\) cells, a cellwise constant coefficient, \(225\) interior
and \(64\) boundary nodes. The port space is the complete nodal trace
space, \(m_\Gamma=64\), so \(q=63\) and the network represents the
complete response space
\(\mathcal V_{i,q}^{\star}(\kappa_i)=\mathcal R_{i,h}(\kappa_i)\). Because
\(\mathcal B_{\Gamma_i}=\Lambda_h(\Gamma_i)\), the exact-block solution is
the fine-grid solution \(u_h\), the reference for every
two-dimensional test. The network has \(18.5\) million parameters and
trained in about \(48\) GPU hours on one NVIDIA H100.

\begin{figure}[!htbp]
\centering
\includegraphics[width=0.8\textwidth]{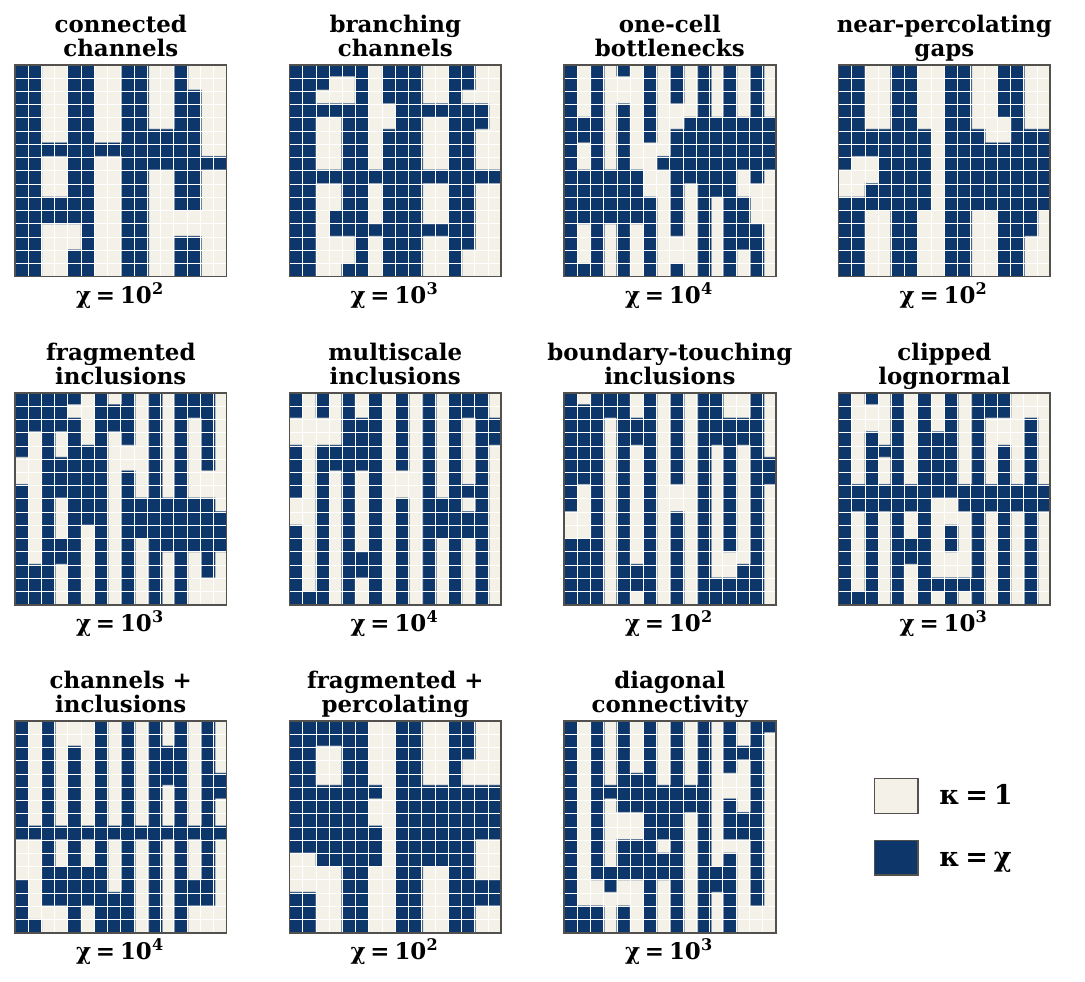}
\caption{One test patch (\(16\times16\) cells) per coefficient
family, in the order of the text; the contrast
\(\chi\in\{10^2,10^3,10^4\}\) cycles across the panels. Cream:
\(\kappa=1\); dark blue: \(\kappa=\chi\). Each panel is the patch of its
family and contrast whose volume fraction of the \(\chi\) phase is closest
to the family median.}
\label{fig:2d-families}
\end{figure}

Every coefficient is binary, \(\kappa\in\{1,\chi\}\) on each cell, with
contrast \(\chi\in\{10^2,10^3,10^4\}\). A coefficient \emph{family}
is a random generator of such layouts with a prescribed connectivity of
the high-conductivity phase; the eleven families are connected channels,
branching channels, one-cell bottlenecks, near-percolating gaps,
fragmented, multiscale and boundary-touching inclusions, clipped
lognormal structures, channels combined with inclusions, fragmented
combined with percolating structures, and diagonally connected layouts; \Cref{fig:2d-families} shows one test patch
of each. Several families contain one-cell-wide channels of
\(\kappa=\chi\) that cross the patch; these channels are what the
classical coarse spaces of the comparison below do not resolve. The training set contains \(1.3\times10^5\) patches;
the validation and test sets contain \(2{,}200\) patches each, drawn from
the same families and contrasts in the same proportions and with the same
distribution of the volume fraction of the \(\chi\) phase. A \emph{mosaic}
is a global domain of \(S\times S\) independently drawn patches:
\emph{mixed} when the patches come from all families and contrasts,
\emph{contrast-restricted} when they share one contrast, \emph{per family}
when they share one family.

\begin{figure}[!htbp]
\centering
\includegraphics[width=\textwidth]{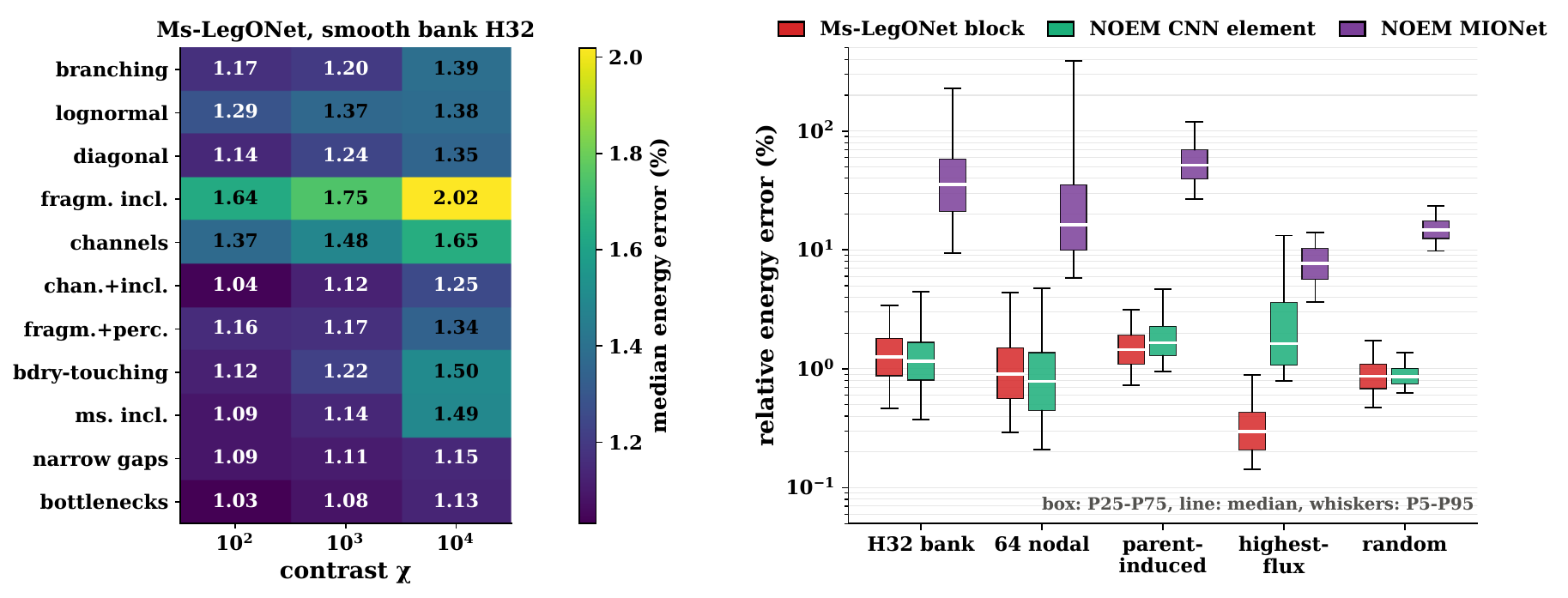}
\caption{Local block quality on the \(2{,}200\) test patches. Left: median
relative energy error of the H32 responses of Ms-LegONet per family
and contrast; the row labels abbreviate the family names of
\Cref{fig:2d-families}: ``channels'' and ``branching'' are the connected
and branching channels, ``narrow gaps'' the near-percolating gaps,
``bottlenecks'' the one-cell bottlenecks, ``fragm.\ incl.'', ``ms.\ incl.''
and ``bdry-touching'' the fragmented, multiscale and boundary-touching
inclusions, ``lognormal'' the clipped lognormal structures,
``chan.+incl.'' and ``fragm.+perc.'' the two combined families, and
``diagonal'' the diagonally connected layouts. Right: distribution (box: P25--P75, line: median,
whiskers: P5--P95) of the response energy error per trace group for
the Ms-LegONet block and for the two NOEM elements of
\Cref{tab:2d-baselines-acc} that read the boundary trace, the
convolutional element and the MIONet of \cite{ouyang2026noem}.}
\label{fig:2d-local}
\end{figure}

\paragraph{Local block and global accuracy}
The learned space is tested first on the \(2{,}200\) test patches
(\Cref{fig:2d-local}) through the relative energy error of the local Galerkin
response \eqref{eq:mslegonet-local-galerkin} against the exact discrete
harmonic extension. Five trace groups are used: \emph{H32}, a fixed bank of
\(32\) piecewise-linear traces; the \(64\) nodal traces;
\emph{parent-induced} traces, which a primitive receives inside \(2\times2\) to
\(8\times8\) blocks of other patches with random external traces; the
four \emph{highest-flux} traces of a patch, which maximize the boundary
flux of the harmonic extension relative to the boundary mass norm;
and \(32\) \emph{random} unit combinations of the port basis per
patch. The
median error is \(0.9\)--\(1.4\%\) for the four generic groups and
\(0.3\%\) for the highest-flux traces, with P95 below \(4.5\%\), and the
H32 median stays between \(1.0\%\) and \(2.1\%\) over all families and
contrasts. The error is the same for traces the block was never shown as
for generic ones, and smallest for the traces that carry the most energy:
the network has learned the response space, not a set of responses.

Global tests use mixed mosaics of \(16\times16\) and \(32\times32\) test
primitives, contrast-restricted mosaics and per-family \(8\times8\) mosaics.
Three problems are solved on every mosaic: \(f=1\) and a Gaussian random
source \(f\), a centred Gaussian random field with covariance operator
\((\tau^2-\Delta)^{-2}\), \(\tau=3\), normalized to unit RMS, both with
\(u=0\) on \(\partial\Omega\), and smooth random Dirichlet data \(g\) with
\(f=0\). Three facts stand out in \Cref{tab:2d-global}. The global energy
error stays at the level of the local one: \(3.4\%\), \(4.1\%\) and \(2.0\%\) for the three
problems on the mixed \(16\times16\) mosaics, within one percentage point of
the \(32\times32\) and per-family \(8\times8\) medians, so composing more
primitives does not accumulate error. The \(L^2\) error is an order of
magnitude smaller, \(0.19\)--\(0.51\%\), because the energy error sits in a
few high-contrast primitives while the solution elsewhere is resolved.
Contrast is the one driver: for the random source the median energy error
rises from \(2.1\%\) at \(\chi=10^2\) and \(2.7\%\) at \(10^3\) to \(8.7\%\)
at \(10^4\) (at \(\chi=10^2,10^3\) every median lies between \(1.6\%\) and
\(2.9\%\), and these rows are omitted from the table).

With the complete nodal port space, \(u_{\mathrm{ms}}^{\star}=u_h\), the first
term of \Cref{prop:localized-bound} vanishes, and \(w_i^{\star}=\mathcal C_{i,h}^{\kappa_i}\gamma_iu_h\) is the
zero-trace correction of the harmonic part of \(u_h\) on \(\omega_i\):
\begin{equation}
e:=\|u_h-u_{\mathrm{ms},\theta}\|_{a_{\kappae}}
\le
\eta_{\mathrm{proj}}
\le
\eta_{\mathrm{ang}}.
\label{eq:2d-cea-bound}
\end{equation}
\Cref{tab:2d-global} reports the effectivity index
\(\eta_{\mathrm{proj}}/e\) per mosaic. It is at least \(1\) on every test,
as \Cref{prop:localized-bound} requires, and its median is
\(1.002\)--\(1.034\) on the mixed layouts and at most \(1.093\) on the
contrast-restricted ones. The global error is therefore, to within a few
percent, the root-sum-square of the local best-approximation errors of the learned
spaces: the skeleton coupling neither amplifies nor cancels them. The angle bound
\(\eta_{\mathrm{ang}}\) is coarser, since it charges every direction of the
learned space with the largest angle. Two implementation checks:
quadtree composition of every depth reproduces the flat solve to
\(3\times10^{-14}\) in relative energy, and the exact-block solution
reproduces \(u_h\) to \(4\times10^{-12}\), both at round-off.

\begin{table}[!htbp]
\centering
\caption[Relative errors of Ms-LegONet and effectivity of the localized bound]{Relative errors (\%) of Ms-LegONet, P50 and P95 over test mosaics
(6 mixed \(16\times16\), 12 mixed \(32\times32\), 3 at \(\chi=10^4\) and
\(16\times16\), 33 per-family \(8\times8\), 6 at \(\chi=10^4\) and
\(32\times32\)), and the effectivity index \(\eta_{\mathrm{proj}}/e\) of
\eqref{eq:2d-cea-bound}, computed per mosaic; min and P50 are taken over the
mosaics. For the layouts with 3 to 12 mosaics the P95 is essentially the
maximum over those mosaics and is to be read as a worst case.}
\label{tab:2d-global}
\footnotesize
\setlength{\tabcolsep}{4pt}
\begin{adjustbox}{max width=\textwidth}
\begin{tabular}{llcccccc}
\toprule
 & & \multicolumn{2}{c}{energy} & \multicolumn{2}{c}{\(L^2\)}
 & \multicolumn{2}{c}{\(\eta_{\mathrm{proj}}/e\)}\\
\cmidrule(lr){3-4}\cmidrule(lr){5-6}\cmidrule(lr){7-8}
layout & problem & P50 & P95 & P50 & P95 & min & P50\\
\midrule
mixed, \(16\times16\) & \(f=1\) & 3.41 & 5.52 & 0.19 & 0.56 & 1.005 & 1.016\\
 & random \(f\) & 4.06 & 4.53 & 0.47 & 0.56 & 1.007 & 1.024\\
 & random \(g\) & 1.95 & 2.60 & 0.51 & 0.74 & 1.001 & 1.002\\
\addlinespace
mixed, \(32\times32\) & \(f=1\) & 4.23 & 5.23 & 0.22 & 0.37 & 1.012 & 1.024\\
 & random \(f\) & 4.59 & 5.55 & 0.47 & 0.75 & 1.020 & 1.034\\
 & random \(g\) & 1.98 & 2.33 & 0.41 & 0.67 & 1.001 & 1.004\\
\addlinespace
\(\chi=10^4\), \(16\times16\) & \(f=1\) & 3.76 & 4.11 & 0.96 & 0.97 & 1.041 & 1.048\\
 & random \(f\) & 8.67 & 10.07 & 3.63 & 6.32 & 1.063 & 1.079\\
 & random \(g\) & 2.02 & 2.14 & 0.64 & 0.71 & 1.001 & 1.001\\
\addlinespace
per family, \(8\times8\) & \(f=1\) & 3.43 & 6.03 & 0.32 & 0.87 & 1.001 & 1.007\\
 & random \(f\) & 3.80 & 7.07 & 0.46 & 1.44 & 1.002 & 1.010\\
 & random \(g\) & 1.97 & 3.05 & 0.46 & 1.03 & 1.000 & 1.002\\
\addlinespace
\(\chi=10^4\), \(32\times32\) & random \(f\) & 7.36 & 7.77 & 3.06 & 3.76 & 1.047 & 1.093\\
\bottomrule
\end{tabular}
\end{adjustbox}
\end{table}

\paragraph{Baselines}
Four baselines are evaluated on a separate set of test benchmark mosaics
(\(16\times16\) test primitives, zero Dirichlet data), every error measured
against the fine-grid solution \(u_h\) on the fine grid. Two are the
classical multiscale methods on the primitive grid \(H=16h\): MsFEM with
linear boundary conditions \cite{hou1997msfem}, and GMsFEM-\(n\) with \(n\)
spectral basis functions per coarse vertex \cite{efendiev2013gmsfem} in the
construction of \cite{rudikov2025gmsfemno}, whose spectral problems are
posed on the four primitives around each vertex and solved anew for every
coefficient. The network of GMsFEM-NO \cite{rudikov2025gmsfemno} learns this
basis from the coefficient, so the GMsFEM-8 rows are the reference its
network is trained to reproduce. Both are local constructions without
oversampling; oversampled and constraint-energy-minimizing variants
\cite{hou1997msfem,chung2023multiscale} remove their boundary-layer
resonance at the price of larger local problems and lie outside the local
setting that GMsFEM-NO and this comparison address. The third is the neural-operator element method (NOEM)
\cite{ouyang2026noem}, the closest design to ours: its element reads
\((\kappa,\text{boundary trace})\) and returns the interior field, and the
elements are composed by minimizing the total energy over their boundary
values with Newton's method. We run it on the same \(7{,}425\) skeleton
unknowns as Ms-LegONet with two elements trained on the same coefficients,
probe traces and exact labels as our block: the MIONet of
\cite{ouyang2026noem}, and the strongest boundary-state element we
obtained on these data, a convolutional element that reads the Coons lift of the trace as an
image and predicts the correction to it with the physical-response loss.
The fourth is an end-to-end FNO \cite{li2021fno} with four Fourier layers,
\(24\) modes and width \(64\), trained on \(10{,}000\) mosaics of training
patches, with the checkpoint selected on held-out mosaics of training
patches, to map the coefficient and the source \(f\) on the \(257^2\) grid to \(u\).

\begin{figure}[!htbp]
\centering
\includegraphics[width=\textwidth]{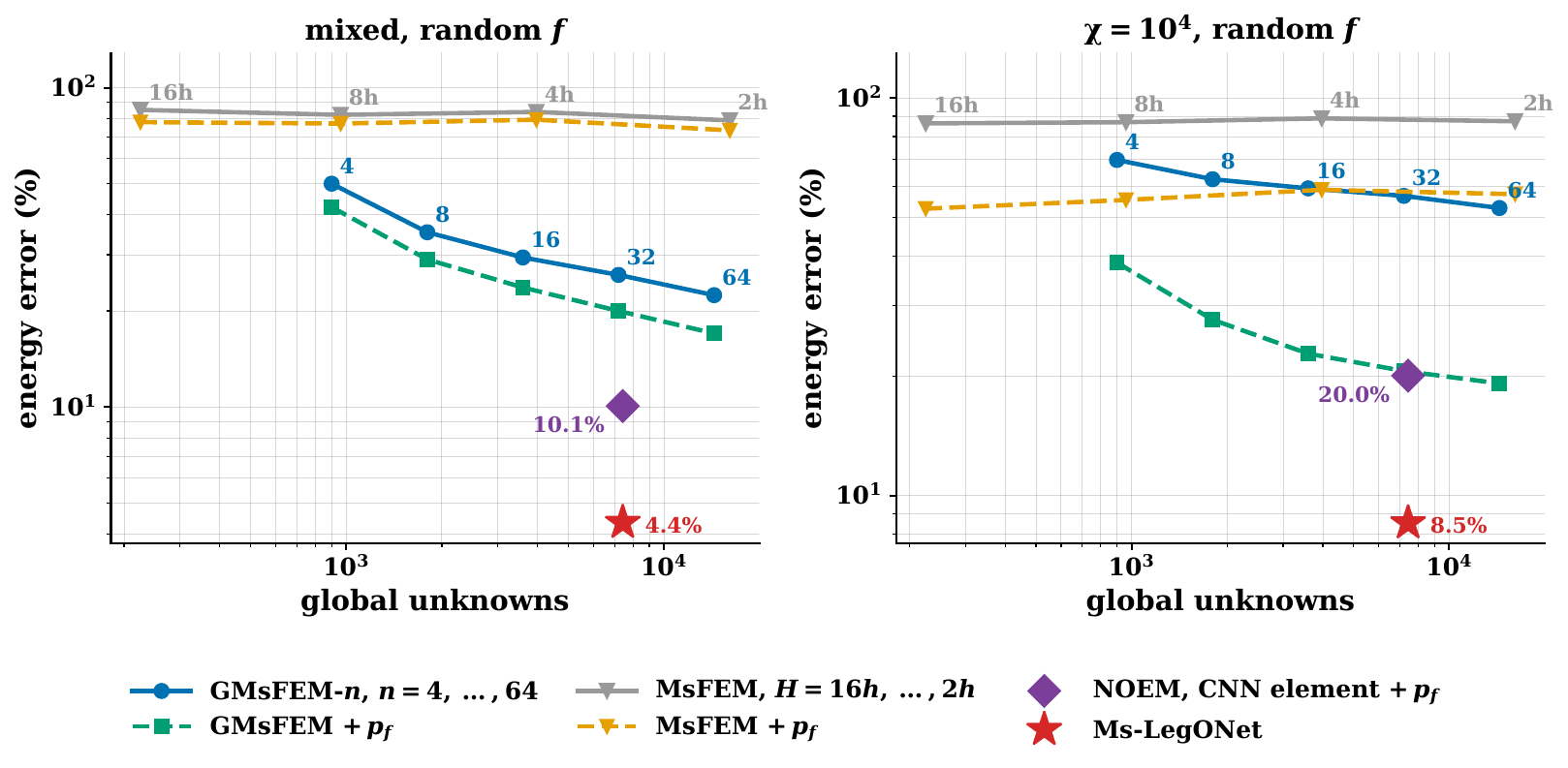}
\caption{Median energy error against the number of global unknowns on the
	mixed and \(\chi=10^4\) test benchmark mosaics with the random source:
	GMsFEM-\(n\) for \(n=4,\ldots,64\), MsFEM under coarse-grid refinement
	\(H=16h,\ldots,2h\), both without and with \(p_f\) (dashed), NOEM with the
	convolutional element \(+\,p_f\), and Ms-LegONet, the last two with
	\(7{,}425\) skeleton unknowns. At \(H=h\) the MsFEM is the fine-grid FEM
	and its error vanishes.}
\label{fig:2d-pareto}
\end{figure}

Two confounders are removed. Each coarse method is also run with the
local particular fields \(p_f\) of Ms-LegONet (``\(+\,p_f\)'': the coarse
Galerkin solution for the load \(\ell_f(v)-a_{\kappae}(p_f,v)\), added to
\(p_f\)), and \Cref{fig:2d-pareto} lets GMsFEM grow through \(n\) and MsFEM
through the coarse-cell size \(H\), past the \(7{,}425\) skeleton unknowns
of Ms-LegONet. What remains is the harmonic space, and four mechanisms
separate the methods (\Cref{tab:2d-baselines-acc}).

\begin{figure}[!htbp]
\centering
\includegraphics[width=\textwidth]{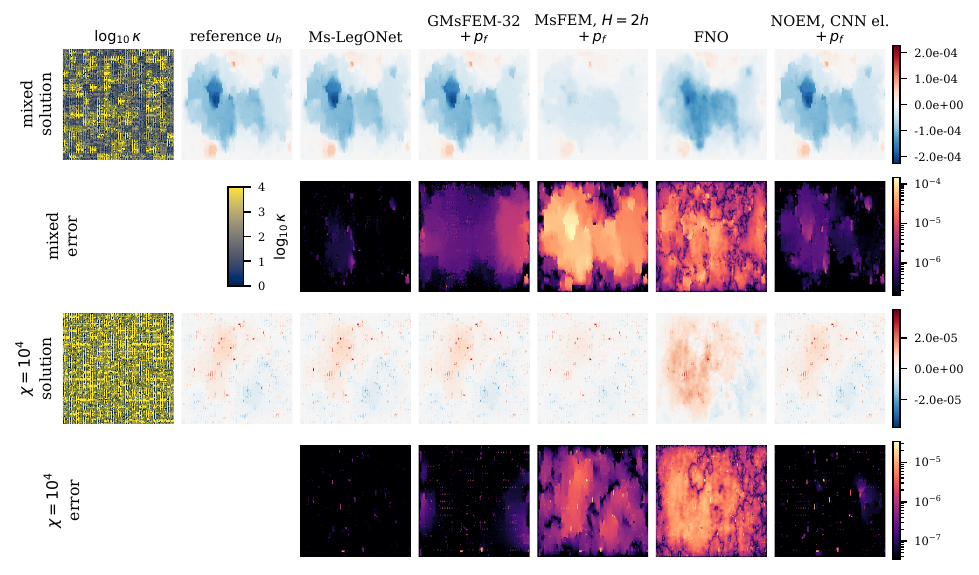}
\caption{Mixed and \(\chi=10^4\) test benchmark mosaics, random source:
	solution (first row per mosaic) and pointwise error \(|u_h-\tilde u|\) on
	a shared logarithmic scale (second row). Columns: \(\log_{10}\kappa\),
	reference \(u_h\), Ms-LegONet, GMsFEM-32 \(+\,p_f\), MsFEM at \(H=2h\)
	\(+\,p_f\), FNO, NOEM with the convolutional element \(+\,p_f\). The
	GMsFEM error is spread over the whole domain, with striations that
	follow the high-conductivity channels of the \(\log_{10}\kappa\) column.}
\label{fig:2d-compare}
\end{figure}

First, in this construction and at these sizes, eigenfunctions per vertex do not resolve one-cell channels. At
equal size, GMsFEM-32 \(+\,p_f\) (\(7{,}200\) unknowns) has \(4.6\) times
the energy error of Ms-LegONet on the mixed mosaics (\(20.0\%\) against
\(4.35\%\)) and \(17\) times the \(L^2\) error (\(8.06\%\) against
\(0.48\%\)); GMsFEM-64 \(+\,p_f\), with twice the unknowns, stays at
\(17\%\). The particular fields matter most at \(\chi=10^4\), where they
take GMsFEM-32 from \(57\%\) to \(21\%\) (\Cref{tab:2d-baselines-acc} and \Cref{fig:2d-pareto}). The
remaining gap, \(21\%\) against the \(8.5\%\) of Ms-LegONet,
we attribute to the channels: in \Cref{fig:2d-compare} the GMsFEM error
follows the high-\(\kappa\) channels across the whole domain.

Second, linear boundary traces do not resolve the coefficient at any \(H>h\).
On the mixed and \(\chi=10^4\) mosaics, from \(H=16h\) to \(H=2h\) with
\(225\) to \(16{,}129\) unknowns, the MsFEM energy error without
\(p_f\) stays between
\(79\%\) and \(89\%\), and its strongest variant, \(H=2h\) with \(p_f\),
at \(57\)--\(74\%\),
still damps the solution to a fraction of its amplitude
(\Cref{fig:2d-compare}). The cause is the boundary condition: wherever a
one-cell channel crosses a coarse edge, the linear trace misrepresents the
solution by an \(O(1)\) amount that does not depend on \(H\), the
boundary-layer mechanism behind the resonance error of
\cite{hou1999convergence} in the periodic setting; the error vanishes
only at \(H=h\), where the method is the fine-grid FEM.

\begin{figure}[!htbp]
\centering
\includegraphics[width=0.56\textwidth]{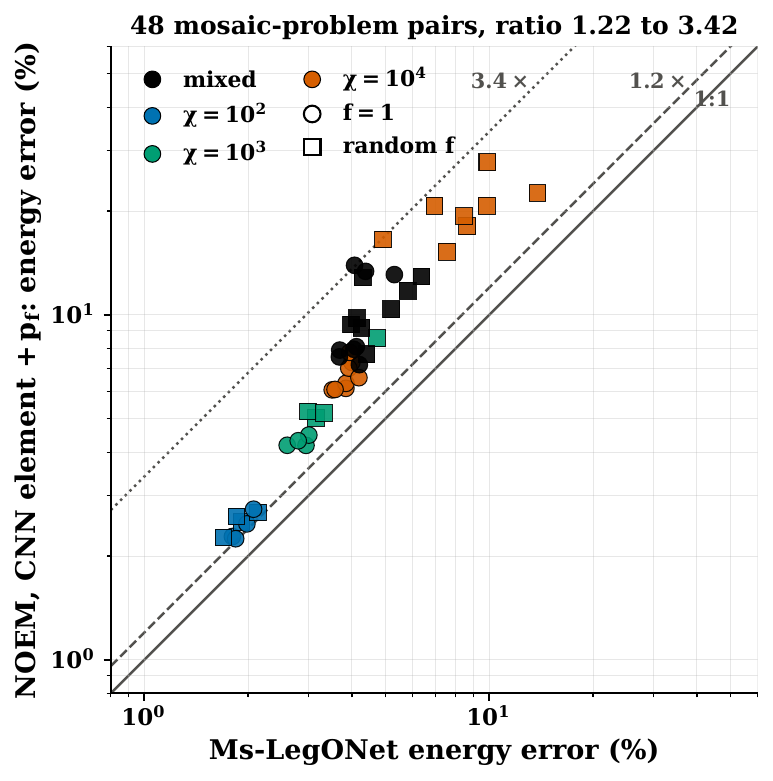}
\caption{Energy error of NOEM with the convolutional element
\(+\,p_f\) against that of Ms-LegONet on the \(48\) mosaic--problem pairs
of the test benchmark (\(24\) mosaics, \(f=1\) and random \(f\)), with the
lines \(y=x\), \(y=1.2x\) and \(y=3.4x\); the ratio lies between \(1.22\)
and \(3.42\) on every pair.}
\label{fig:2d-noem-scatter}
\end{figure}

Third, the elements that read the trace are accurate only on the trace
distribution they saw. Locally, on the \(2{,}200\) test patches of
\Cref{fig:2d-local} (right panel), the MIONet element has \(14.7\%\) median
energy error on random traces and \(51\%\) on the parent-induced traces of
a global solve. The convolutional element matches our block on the trace
groups it was trained on (\(0.86\%\) against \(0.87\%\) on random
traces, \(0.79\%\) against \(0.91\%\) on the nodal traces) and loses it on
the highest-flux traces it never saw (\(1.63\%\), P95 \(13\%\), against
\(0.30\%\) and \(0.89\%\)). Globally (\Cref{tab:2d-baselines-acc}), the
MIONet element gives \(54\)--\(71\%\) energy error, and the convolutional
element \(10.1\%\) on the mixed mosaics and \(20\%\) at \(\chi=10^4\) with
the random source. In \Cref{fig:2d-pareto} this point lies at or below the GMsFEM
\(+\,p_f\) curve; mosaic by mosaic (\Cref{fig:2d-noem-scatter}), its
energy error is \(1.2\) to \(3.4\) times that of Ms-LegONet on all
\(48\) mosaic--problem pairs, with the smallest ratios at \(\chi=10^2\). The gap is the element, not the optimizer: Newton started
at the exact skeleton values of \(u_h\) ends at the same energy on \(46\)
of \(48\) solves, and fed the true traces the element's response has
\(12.5\%\) error against \(4.9\%\) for the composed Ms-LegONet responses
(the term \(\eta_{\mathrm{proj}}\) of \Cref{prop:localized-bound}).
In \Cref{fig:2d-compare} the inexact traces propagate through the
energy minimization and spread the NOEM error over the whole mixed
mosaic, whereas the Ms-LegONet error is confined to a few primitives: the
block reads \(\kappa\) alone, and the trace enters only through the local
Galerkin solve \eqref{eq:mslegonet-local-galerkin}, which introduces no
approximation error in the trace. The
composition also differs in kind: the NOEM energy Hessian was indefinite in
three quarters of the solves, which took a median of \(16\) to \(24\)
Newton steps with a skeleton factorization each, against one SPD solve
for Ms-LegONet, and the nonlinear composition admits no bound of the form
of \Cref{prop:localized-bound}.

Fourth, the global map offers the FNO little to exploit. The FNO reaches
\(38\%\) median \(L^2\) error on the mixed mosaics and \(457\%\) in energy
(\Cref{tab:2d-baselines-acc}); the solution of \(256\) independent patches
with one-cell channels has little low-frequency structure that \(24\) Fourier
modes can represent, and in \Cref{fig:2d-compare} its error has the size
of the solution everywhere.

Across the four baselines, in \Cref{fig:2d-compare} the Ms-LegONet error lies an
order of magnitude below GMsFEM-32 and two below MsFEM over most of the
domain and is confined to a few high-contrast primitives. With
\(f=1\) (the \(f=1\) column of \Cref{tab:2d-baselines-acc}) the ordering is
the same: \(4.1\%\) energy error against \(8.0\%\) for NOEM with the
convolutional element, \(23\%\) for GMsFEM-32 \(+\,p_f\) and \(76\%\) for
the strongest MsFEM variant.

\begin{table}[!htbp]
	\centering
	\caption{Median relative energy / \(L^2\) error (\%) against the fine-grid
		solution \(u_h\) on the test benchmark mosaics (\(16\times16\) primitives;
		8 mixed, 8 at \(\chi=10^4\), 4 at \(\chi=10^2\), 4 at \(\chi=10^3\)),
		random source \(f\) (all columns) and \(f=1\) (the \(f=1\) column,
		mixed mosaics). ``unknowns'': size of the skeleton or coarse system. Each
		classical method appears in its basic form on the primitive grid
		\(H=16h\) and in its strongest form, with the most unknowns and the
		particular fields \(p_f\); GMsFEM-32 is listed without and with
		\(p_f\). NOEM: Newton energy minimization over the
		same skeleton values with an element that reads \((\kappa,\) boundary
		trace\()\), with \(p_f\); the MIONet row uses the architecture and loss
		of \cite{ouyang2026noem}, the CNN row our strongest boundary-state
		element.}
	\label{tab:2d-baselines-acc}
	\footnotesize
	\setlength{\tabcolsep}{3.2pt}
	\begin{adjustbox}{max width=\textwidth}
		\begin{tabular}{lcccccc}
			\toprule
			 & & \multicolumn{2}{c}{mixed} & & & \\
			\cmidrule(lr){3-4}
			method & unknowns & random \(f\) & \(f=1\) & \(\chi=10^2\) & \(\chi=10^3\) & \(\chi=10^4\)\\
			\midrule
			Ms-LegONet & 7{,}425 & \textbf{4.35 / 0.48} & \textbf{4.10 / 0.28} & \textbf{1.88 / 0.11} & \textbf{3.23 / 0.43} & \textbf{8.54 / 5.78}\\
			MsFEM & 225 & 85.21 / 71.21 & 80.32 / 65.55 & 58.00 / 28.79 & 69.42 / 45.66 & 86.31 / 82.64\\
			MsFEM, \(H=2h\), \(+\,p_f\) & 16{,}129 & 73.53 / 65.26 & 76.27 / 58.88 & 59.30 / 40.08 & 68.54 / 55.34 & 57.29 / 53.24\\
			GMsFEM-8 & 1{,}800 & 35.28 / 12.58 & 30.90 / 9.92 & 24.43 / 5.64 & 35.93 / 12.57 & 62.47 / 50.37\\
			GMsFEM-32 & 7{,}200 & 25.93 / 8.93 & 23.69 / 7.12 & 16.36 / 3.49 & 27.83 / 9.02 & 56.72 / 46.98\\
			GMsFEM-32 \(+\,p_f\) & 7{,}200 & 20.00 / 8.06 & 23.36 / 7.01 & 14.38 / 3.22 & 19.03 / 6.29 & 20.58 / 16.15\\
			GMsFEM-64 \(+\,p_f\) & 14{,}400 & 17.01 / 6.87 & 20.46 / 5.69 & 11.75 / 2.53 & 16.12 / 5.00 & 19.15 / 15.61\\
			NOEM, MIONet \(+\,p_f\) & 7{,}425 & 68.42 / 54.88 & 70.71 / 50.88 & 54.12 / 32.76 & 65.16 / 48.36 & 55.83 / 51.28\\
			NOEM, CNN element \(+\,p_f\) & 7{,}425 & 10.08 / 2.12 & 8.02 / 0.87 & 2.56 / 0.21 & 5.22 / 0.96 & 20.04 / 20.15\\
			FNO & \(257^2\) grid & 457.35 / 38.12 & 335.43 / 20.90 & 101.01 / 59.68 & 218.09 / 36.19 & 736.07 / 130.39\\
			\bottomrule
		\end{tabular}
	\end{adjustbox}
\end{table}

\paragraph{Many coefficients, one quantity of interest}
\Cref{tab:2d-homog} reports a many-query study on \(20\) mixed test mosaics per
layout size \(S\in\{8,16,32,64\}\). For each coefficient the effective
conductivity tensor \(K^\ast_{ij}=|\Omega|^{-1}a_\kappa(u^i,u^j)\) is
computed from the two cell problems \(f=0\), \(u^j=x_j\) on
\(\partial\Omega\), with \(x_j\) the \(j\)-th coordinate function: the
nonhomogeneous Dirichlet problem of \Cref{subsec:source-global} without
source. The reference is the monolithic \(Q_1\) solve by PARDISO; AMG-PCG
agrees with it to below \(10^{-4}\%\) in \(K^\ast\) and is listed for the
cost of an iterative solver at this contrast. At \(S=64\) (\(4{,}096\)
primitives) Ms-LegONet needs \(2.3\) s per coefficient against \(3.8\) s
for PARDISO and \(57\) s for AMG-PCG, and its \(K^\ast\) error,
\(\max_{ij}|K^\ast_{ij}-K^\ast_{ij,h}|/\|K^\ast_h\|_2\), is
\(0.08\)--\(0.19\%\) in median and at most \(0.26\%\), although the cell
solutions carry a \(4.45\%\) energy error. The two numbers are consistent:
\(K^\ast_{jj}\) is the energy of a Galerkin solution, so by Galerkin
orthogonality its relative error is the squared relative energy error,
\((4.45\%)^2\approx0.2\%\), and \(K^\ast_{ij}\) follows by polarization.

\begin{table}[!htbp]
\centering
\caption[Many-query effective conductivity]{Many-query effective conductivity on test mosaics: time per
coefficient (s; median over 20 queries) and \(K^\ast\) error of
Ms-LegONet against FEM with PARDISO (median / max over the queries, \%).
AMG-PCG agrees with PARDISO to below \(10^{-4}\%\) in \(K^\ast\).}
\label{tab:2d-homog}
\footnotesize
\setlength{\tabcolsep}{5pt}
{
\begin{tabular}{lcccc}
\toprule
 & \multicolumn{4}{c}{\(S\)}\\
\cmidrule(lr){2-5}
 & \(8\) & \(16\) & \(32\) & \(64\)\\
\midrule
\multicolumn{5}{l}{\emph{time per coefficient (s)}}\\
Ms-LegONet     & \textbf{0.038} & \textbf{0.137} & \textbf{0.569} & \textbf{2.287}\\
FEM, PARDISO   & 0.086 & 0.219 & 0.936 & 3.809\\
FEM, AMG-PCG   & 0.745 & 3.325 & 14.184 & 57.267\\
\addlinespace
\multicolumn{5}{l}{\emph{\(K^\ast\) error of Ms-LegONet (\%), median / max}}\\
               & 0.08 / 0.19 & 0.15 / 0.26 & 0.16 / 0.25 & 0.19 / 0.24\\
\bottomrule
\end{tabular}
}
\end{table}
\subsection{Three-dimensional layered media with reduced ports}
\label{subsec:3d-layered}

The three-dimensional study uses reduced boundary ports, so both terms of
\eqref{eq:num-error-split} are present, and its coefficient class admits an
exact continuum reference.

\paragraph{Problem and primitive block}
The domain is the \emph{chain} \(\Omega_N=(0,N)\times(0,1)^2\),
partitioned into \(N\) unit primitives \(\omega_i=(i-1,i)\times(0,1)^2\), with
homogeneous Dirichlet data on \(\partial\Omega_N\). The sources are
\(f=1\), which excites only the transverse modes with \(m\) and \(n\) odd
(four of them lie in the port space below), and
\(f_2=(1+y)(1+2z)\), which has no symmetry in \(y\) or \(z\) and excites
every mode. On each primitive \(\kappa\) is
piecewise constant on \(16^3\) cells and varies in \(x\), in \(16\)
layers with values drawn independently from the log-uniform ladder
\(\{\chi^{j/J}:j=0,\ldots,J\}\), \(J\in\{2,3,4,6,8\}\) chosen per patch,
so \(\kappa\) ranges over \([1,\chi]\) with contrast
\(\chi\in\{10^2,10^3,10^4\}\). Primitives draw their layer sequences
independently, so \(\kappa\) is discontinuous across interfaces; a
\emph{mixed} chain draws from all three contrasts, a contrast-restricted chain
from one. The port space on each \(x\)-face is
\begin{equation}
\mathcal B_{\Gamma}^{x}
=\operatorname{span}\{\sin(m\pi y)\sin(n\pi z):\ 1\le m,n\le 4\},
\label{eq:3d-port-space}
\end{equation}
so \(m_\Gamma=32\) ports per primitive; the four lateral faces lie on
\(\partial\Omega_N\) and carry no ports. The lifting is \(E_ig=(1-\hat x)g\)
or \(\hat x g\) by face, with \(\hat x=x-(i-1)\) the local coordinate; interface traces outside
\eqref{eq:3d-port-space} are set to zero.
The learned block maps the normalized \(\log\kappa\) on the \(16^3\) cells to a
\(q=32\) dimensional zero-trace space. Each predicted function is a sum
of \(4\) separated terms,
\(\widetilde\psi_{\theta,r}(x,y,z)=\sum_{p=1}^{4}
X_{rp}(x)\,Y_{rp}(y)\,Z_{rp}(z)\),
whose one-dimensional factors are continuous piecewise quintics with zero
end values; every term therefore lies in \(H^1_0\), and the energy Gram
matrices reduce to products of one-dimensional integrals, evaluated
exactly. The network has \(3.7\)
million parameters and trains in about \(4.6\) GPU hours.
Because \(\kappa=\kappa(x)\), the sine modes diagonalize the operator in
\((y,z)\): each transverse mode \((m,n)\) reduces to a one-dimensional
problem on \((0,N)\), solved
exactly layer by layer with exponential elements. The reference is this
continuum solution. The labels, the \(32\) continuum corrections
\(\mathcal C_i^{\kappa_i}b_\alpha^{(i)}\) of a primitive (the continuum
analogue of \eqref{eq:local-correction-operator}), are computed offline
with the exact modal solver, and so are the local particular fields
\(p_i(f)\) of \eqref{eq:local-particular-problem}, which contain no learned
quantity (\Cref{app:q1-pf} reports the effect of computing \(p_f\) on the
\(Q_1\) mesh instead). The separability enters these local computations
and the one-dimensional energy integrals of the separated representation;
the network sees the raw \(16^3\) array of \(\kappa\), the learned block
replaces the harmonic part alone, and the assembly is the general one of
\Cref{subsec:composition}. The training,
validation and test sets contain \(5{,}500\), \(550\) and \(550\) patches.

\begin{figure}[!htbp]
\centering
\includegraphics[width=0.9\textwidth]{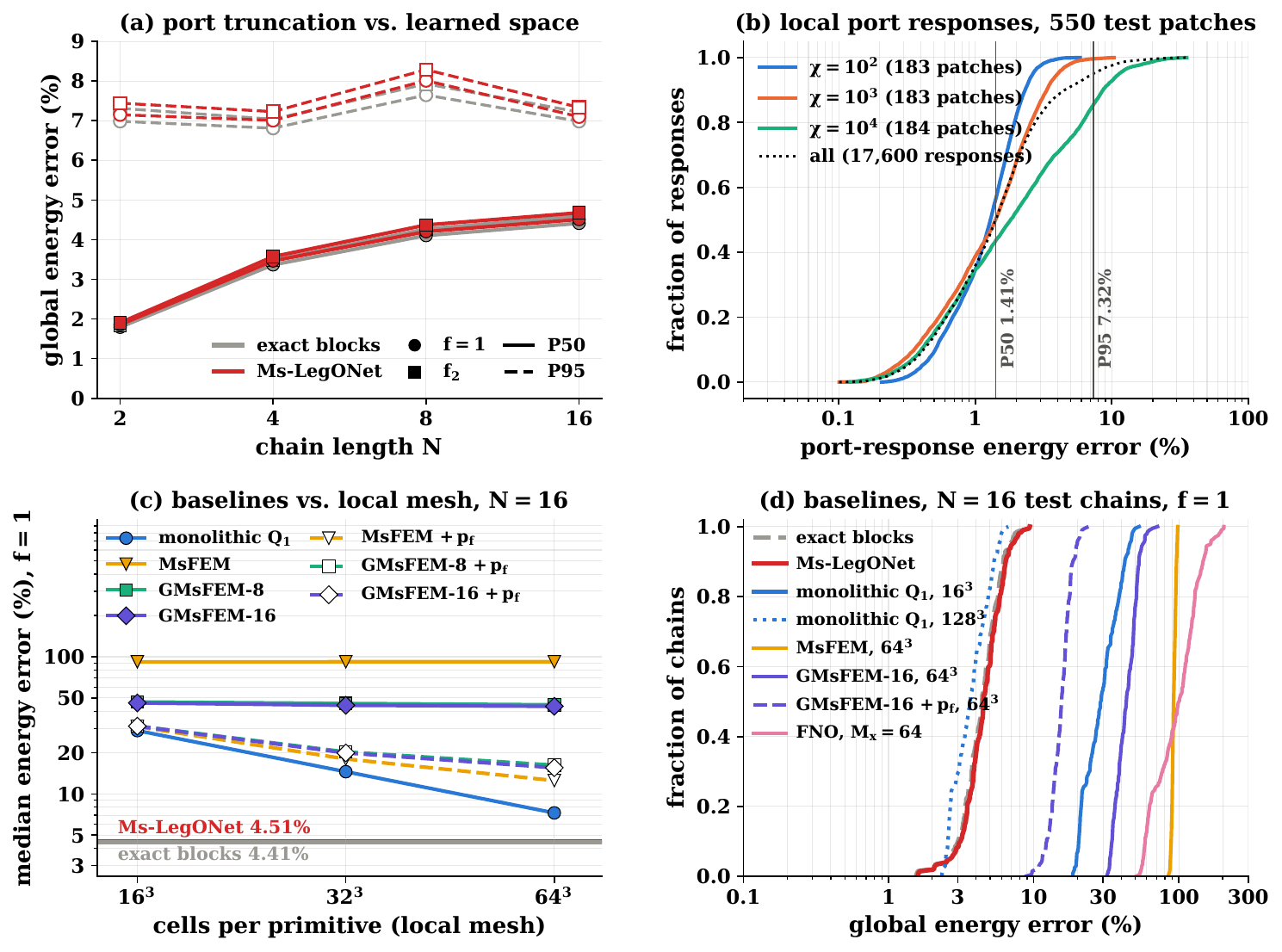}
\caption{Accuracy of the three-dimensional block. (a) Global energy
error against the continuum reference on the \(800\) test chains (\(200\)
per \(N\)): exact blocks and Ms-LegONet, P50 (solid) and P95 (dashed),
\(f=1\) (circles) and \(f_2\) (squares); the two medians sit \(0.1\) percentage
points apart at every \(N\). (b) Empirical distribution of the energy
error of the \(550\times32\) test port responses by contrast; vertical
lines: overall P50 (\(1.41\%\)) and P95 (\(7.32\%\)). (c) Median energy
error (\(f=1\)) of the monolithic \(Q_1\) solution and of MsFEM and
GMsFEM-\(n\), without and with \(p_f\) (dashed), against the local mesh of
\(16^3\), \(32^3\) and \(64^3\) cells per primitive on the \(200\) test
chains with \(N=16\); horizontal lines: Ms-LegONet (\(4.51\%\)) and the
exact blocks (\(4.41\%\)). (d) Empirical distribution of the global
energy error over the same \(200\) chains (\(f=1\)): exact blocks,
Ms-LegONet, the monolithic \(Q_1\) solution on \(16^3\) and \(128^3\)
cells per primitive, MsFEM, GMsFEM-16 and GMsFEM-16 \(+\,p_f\) on the
\(64^3\) local mesh, and the FNO with \(M_x=64\).}
\label{fig:3d-accuracy}
\end{figure}

\paragraph{Accuracy}
Global tests use \(800\) chains assembled from the \(550\) test patches:
\(N\in\{2,4,8,16\}\), each with \(50\) mixed chains and \(50\) at each of
\(\chi=10^2,10^3,10^4\). The two terms of \eqref{eq:num-error-split} behave
differently (\Cref{fig:3d-accuracy}a). The port truncation sets a
floor that saturates with \(N\): with \(f=1\) the exact-block median
grows from \(1.8\%\) at \(N=2\) to \(4.1\%\) at \(N=8\) and \(4.4\%\) at
\(N=16\), while its 95th percentile stays between \(6.8\%\) and
\(7.7\%\); at \(N=16\) the two sources give medians of
\(4.4\)--\(4.6\%\) and 95th percentiles of \(7.0\)--\(7.2\%\)
(\Cref{tab:3d-baselines}). The
discarded modes, \(\max(m,n)\ge5\),
decay from each interface as \(e^{-\xi_{mn}|x-x_\Gamma|}\) with
\(\xi_{mn}:=\pi\sqrt{m^2+n^2}\ge\pi\sqrt{26}\approx16\), so the truncation
error lives in thin layers at the interfaces and does not propagate along
the chain. The learned space adds \(0.1\) percentage points to the
median of this
floor at every \(N\) (\Cref{tab:3d-baselines}; the median curves in \Cref{fig:3d-accuracy}a
and the distributions in \Cref{fig:3d-accuracy}d coincide; at the 95th percentile the gap reaches
\(0.4\) points at \(N=8\)). On the \(550\times32\) test port responses its
energy error has a median of \(1.41\%\) and a 95th percentile of
\(7.32\%\), the tail from the \(\chi=10^4\) patches (\Cref{fig:3d-accuracy}b:
P95 of \(2.6\%\), \(4.0\%\) and \(11.3\%\) at \(\chi=10^2\), \(10^3\) and
\(10^4\)). Globally, Ms-LegONet
reaches \(4.51\%\) against \(4.41\%\) for the exact blocks at \(N=16\) with
\(f=1\), \(7.10\%\) against \(6.98\%\) at the 95th percentile, and
\(4.68\%\) against \(4.60\%\) with \(f_2\); over all \(800\) chains,
\(N=2\) to \(16\), the medians are \(3.83\%\) for Ms-LegONet and \(3.73\%\) for the exact
blocks.
The two numbers are consistent: the learned-space error is nearly
orthogonal to the port truncation in energy, so the two add approximately
in quadrature, and \(\sqrt{4.41^2+e^2}=4.51\) corresponds to a learned-space
error \(e\) of about \(1\%\), close to the \(1.4\%\) median port-response
error of the block.
The accuracy of Ms-LegONet is that of the port space, which is chosen, not
learned. Global errors are evaluated exactly from the energy of the
continuum reference through
\(\|u-u_{\mathrm{ms},\theta}\|_{a_{\kappae}}^2
=a_{\kappae}(u,u)-2\ell_f(u_{\mathrm{ms},\theta})
+a_{\kappae}(u_{\mathrm{ms},\theta},u_{\mathrm{ms},\theta})\).
Composing an \(N=16\) chain pairwise, \(2\to4\to8\to16\), reproduces the
single-level solve to \(3.9\times10^{-12}\) on \(40\) test chains.

\begin{figure}[!htbp]
\centering
\includegraphics[width=\textwidth]{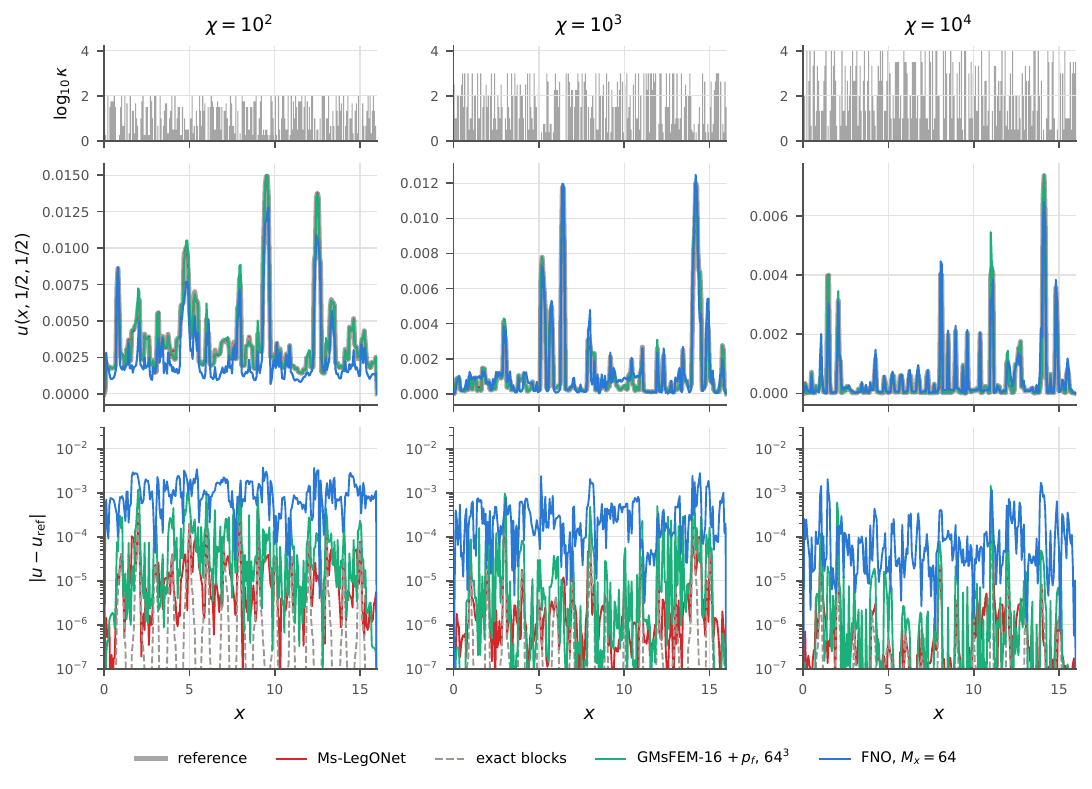}
\caption{Three test chains with \(N=16\), one per contrast group
(columns: \(\chi=10^2\), \(10^3\), \(10^4\); \(f=1\); the chain with the
median Ms-LegONet error of its group). Top: \(\log_{10}\kappa\) of the
\(256\) layers. Middle: solution on the centre line \(y=z=1/2\): continuum
reference, exact blocks, Ms-LegONet, GMsFEM-16 \(+p_f\) on the \(64^3\) local
mesh and the FNO with \(M_x=64\). Bottom: pointwise error against the
continuum reference.}
\label{fig:3d-compare}
\end{figure}

\paragraph{Baselines}
\Cref{tab:3d-baselines}, \Cref{fig:3d-compare} and
\Cref{fig:3d-accuracy}c,d
compare Ms-LegONet on the \(200\) test chains with \(N=16\) against five
baselines. The \emph{\(Q_1\) blocks} are the Ms-LegONet construction with
\(Q_1\) local solves for both the port responses and \(p_f\) on the \(16^3\)-cell primitive mesh in place of the
network. The monolithic \(Q_1\) solution is refined from \(16^3\) to
\(128^3\) cells per primitive and solved by algebraic multigrid (smoothed
aggregation, PCG to a relative residual of \(10^{-8}\)). MsFEM uses one
\(\kappa\)-harmonic basis function per interface, the extension of the
lowest port trace \(\sin\pi y\sin\pi z\). GMsFEM-\(n\) uses \(n\)
eigenfunctions per interface of the local spectral problem on the two
adjacent primitives (Dirichlet on the lateral faces, Neumann on the
\(x\)-faces), times the one-dimensional \(\kappa\)-harmonic hat as partition
of unity; GMsFEM-16 has the skeleton size of Ms-LegONet. Both are run with
and without \(p_f\) on local meshes of \(16^3\), \(32^3\) and \(64^3\)
cells, exact to round-off by separability. The FNO has four spectral layers
of width \(32\) with \((M_x,8,8)\) Fourier modes, \(M_x\in\{16,32,64\}\)
(\(34\)--\(134\)M parameters), reads \((\log\kappa,f,x,y,z)\) on the
\(257\times17\times17\) nodal grid of an \(N=16\) chain, and is trained on
\(10{,}000\) chains of training patches with continuum labels and both
sources. No classical or monolithic baseline is restricted to the
\(16^3\) mesh on which the block reads \(\kappa\).

Three facts emerge. First, the monolithic error is first order in
\(h\), and its price grows accordingly: \(29\%\) at \(0.5\) s, \(15\%\) at \(3.7\) s,
\(7.3\%\) at \(35\) s and \(3.6\%\) at \(406\) s per chain at \(N=16\) with
\(f=1\) (\Cref{tab:3d-baselines}; \Cref{fig:3d-accuracy}d), against \(4.5\%\)
at \(0.028\) s for Ms-LegONet. Matching the
learned block takes about \(100^3\) cells per primitive and three to four
orders of magnitude more time. Both the learned and the exact blocks
use the exact \(p_f\) here; with \(p_f\) computed on the \(Q_1\) mesh
instead, the global error of Ms-LegONet is \(29.6\%\), \(15.4\%\) and
\(8.8\%\) on \(16^3\), \(32^3\) and \(64^3\) cells, the \(Q_1\) error of
\(p_f\) added in quadrature to the \(4.4\%\) port floor
(\Cref{app:q1-pf}): the learned harmonic space is not the bottleneck. Second, the coarse spaces do not improve
with the local mesh (\Cref{fig:3d-accuracy}c): from \(16^3\) to \(64^3\) cells MsFEM stays at \(92\%\)
and GMsFEM-16 moves from \(46\%\) to \(44\%\) (\(47\%\) to \(45\%\) for
GMsFEM-8); the particular fields of the refined mesh take them from
\(31\%\) to \(12.5\%\) and \(15.5\%\) (\(64^3\) rows of
\Cref{tab:3d-baselines}), three times the error of Ms-LegONet
at the same skeleton size. Neither does the number of modes help: GMsFEM-8
and GMsFEM-16 coincide to within two points for \(f=1\), and with \(p_f\) the
single-mode MsFEM ends below GMsFEM-16 (\(12.5\%\) against \(15.5\%\)),
so the two spaces are not nested. They fail for different reasons. The
single MsFEM function is the exact \(\kappa\)-harmonic response of the
lowest port mode, so with \(p_f\) its error is the truncation to one mode:
\(12.5\%\) on \(64^3\) cells is consistent with the \(7.3\%\) \(Q_1\)
error of \(p_f\) and a one-mode truncation of about \(10\%\) in
quadrature. The GMsFEM functions have their \(x\)-dependence fixed by the
partition of unity and the spectral problem and lack the exponential decay
\(e^{-\xi_{mn}x}\) of a harmonic extension, so sixteen of them resolve
less than one exact response. The learned space, a free
\(32\)-dimensional zero-trace space, has neither constraint.
Third, the FNO, with \(9\) to \(36\) times the
parameters of the local block, reaches \(15.7\%\), \(13.2\%\) and \(11.5\%\)
median nodal \(L^2\) error on the hold-out chains for \(M_x=16,32,64\), and
\(34\%\), \(31\%\) and \(28\%\) on the test chains, against
\(2.2\%\) for Ms-LegONet in the same norm (last column of
\Cref{tab:3d-baselines}). In this norm
the \(Q_1\) discretization error on the \(16^3\) mesh is
\(1.5\)--\(1.6\%\), against \(29\%\) in energy: the \(L^2\) error of a
\(Q_1\) solution converges at second order and its energy error at first
(a factor of four against two per mesh halving in \Cref{tab:3d-baselines}),
so the nodal \(L^2\) norm flatters every \(Q_1\) row and the energy columns remain the measure of the field. The FNO output, read as a
\(Q_1\) field on that mesh, has \(100\%\) energy error, three times the
\(29\%\) of the \(Q_1\) solution (\Cref{tab:3d-baselines}). On the
centre line \(y=z=1/2\) of three \(N=16\) test chains
(\Cref{fig:3d-compare}) the pointwise error of Ms-LegONet is
\(10^{-6}\)--\(10^{-4}\), matches the exact blocks at the interfaces and
lies one to two orders of magnitude below the FNO error, which is visible
in the peaks of the low-conductivity layers.

\begin{table}[!htbp]
\centering
\caption{Baselines on the 200 test chains with $N=16$: relative energy error (\%) against the continuum reference, P50 and P95 for $f=1$ and $f_2$; wall time per solve for a new chain of $N=16$ primitives ($f=1$, median of five runs on the node of \Cref{sec:numerical-results}); and the relative nodal $L^2$ error (\%, P50, $f=1$) on the $257\times17\times17$ nodal grid of the $16^3$ mesh, on which the $Q_1$ error is of second order. Exact blocks: continuum port responses from the modal solver, specific to layered media and not timed. $Q_1$ blocks: the Ms-LegONet construction with $Q_1$ local solves for both the port responses and $p_f$ on $16^3$ cells per primitive. Monolithic $Q_1$: algebraic multigrid on $h=1/16$ to $1/128$. MsFEM and GMsFEM-$n$ on the $64^3$ local mesh, evaluated through the separable structure of $\kappa$, which is specific to layered media and is not timed; a generic implementation on the $16^3$ mesh takes $26.0$ s for GMsFEM-8 and $1.35$ s for MsFEM. FNO with $(M_x,8,8)$ modes: output read as a $Q_1$ field on the $16^3$ mesh, which no $Q_1$ field can bring below the $28.9\%$ of the monolithic $h=1/16$ row.}
\label{tab:3d-baselines}
\footnotesize
\begin{adjustbox}{max width=\textwidth}
\begin{tabular}{l rr rr r r}
\toprule
& \multicolumn{2}{c}{$f=1$ [\%]} & \multicolumn{2}{c}{$f_2$ [\%]} & time [s] & nodal $L^2$ [\%]\\
\cmidrule(lr){2-3}\cmidrule(lr){4-5}\cmidrule(lr){6-6}\cmidrule(lr){7-7}
method & P50 & P95 & P50 & P95 & $N=16$ & P50\\
\midrule
\textbf{Ms-LegONet} & \textbf{4.51} & \textbf{7.10} & \textbf{4.68} & \textbf{7.33} & \textbf{0.028} & \textbf{2.19} \\
exact blocks & 4.41 & 6.98 & 4.60 & 7.22 & -- & 2.16 \\
$Q_1$ blocks, $16^3$ & 29.44 & 45.91 & 30.08 & 46.69 & 1.61 & 1.56 \\
monolithic $Q_1$, $16^3$ ($h=1/16$) & 28.94 & 45.78 & 29.58 & 46.56 & 0.52 & 1.50 \\
monolithic $Q_1$, $32^3$ ($h=1/32$) & 14.52 & 22.97 & 14.85 & 23.37 & 3.67 & 0.34 \\
monolithic $Q_1$, $64^3$ ($h=1/64$) & 7.27 & 11.50 & 7.44 & 11.71 & 35.3 & 0.08 \\
monolithic $Q_1$, $128^3$ ($h=1/128$) & 3.64 & 5.75 & 3.72 & 5.86 & 406.2 & 0.02 \\
\addlinespace
MsFEM, $64^3$ & 91.92 & 97.07 & 92.42 & 97.25 & -- & 88.70 \\
GMsFEM-8, $64^3$ & 44.63 & 58.78 & 49.61 & 62.92 & -- & 26.85 \\
GMsFEM-16, $64^3$ & 43.52 & 56.39 & 46.30 & 59.20 & -- & 26.08 \\
MsFEM $+\,p_f$, $64^3$ & 12.51 & 17.71 & 15.08 & 20.50 & -- & 6.15 \\
GMsFEM-8 $+\,p_f$, $64^3$ & 16.20 & 21.06 & 17.98 & 22.75 & -- & 7.01 \\
GMsFEM-16 $+\,p_f$, $64^3$ & 15.53 & 19.93 & 16.33 & 21.19 & -- & 6.86 \\
\addlinespace
FNO, $M_x=16$ & 164.07 & 225.58 & 152.70 & 218.64 & 0.005 & 34.31 \\
FNO, $M_x=32$ & 117.39 & 183.42 & 117.86 & 163.24 & 0.004 & 30.60 \\
FNO, $M_x=64$ & 99.96 & 165.22 & 102.87 & 158.02 & 0.004 & 28.31 \\
\bottomrule
\end{tabular}
\end{adjustbox}
\end{table}

\paragraph{Many coefficients: the effective port response of a chain}
We time the \(32\times32\) end-to-end response
\(\mathbf D\) of a chain, \(\mathbf D_{\alpha\beta}=a_\kappa(u_{\alpha},u_{\beta})\), where \(u_{\alpha}\)
is the \(\kappa\)-harmonic field with port trace \(b_{\alpha}\) on one end face and
zero on the other, for \(50\) new chains per \(N\in\{4,8,16,32,64\}\).
\Cref{tab:3d-eff} compares Ms-LegONet (network, Grams, Schur blocks,
condensation) with the \(Q_1\) blocks on \(16^3\) cells per primitive, all
primitives at once (assembly, Cholesky factorization of the interior
blocks, port solves, Schur blocks), and with the monolithic \(Q_1\) solve on
\(16^3\) and on \(32^3\) cells per primitive. The monolithic solve builds
one algebraic multigrid hierarchy of the chain and performs \(32\) end-face
solves by preconditioned CG to a relative residual of \(10^{-8}\). The exact blocks, whose modal
responses exist only for layered media, are not timed; in a general medium
the exact block is the \(Q_1\) block. For \(64\) primitives
Ms-LegONet computes \(\mathbf D\) in \(0.11\) s against \(6.4\) s
for the \(Q_1\) blocks, \(11.5\) s for the monolithic solve on the same
mesh and \(98.5\) s on \(32^3\) cells: \(58\), \(106\) and \(904\) times
faster. Over \(N=4\) to \(64\) the factors are \(31\)--\(60\), \(59\)--\(115\)
and \(388\)--\(904\).

The learned \(\mathbf D\) is also more accurate, by a factor of about
\(600\). Its relative error has median \(7\times10^{-5}\)--\(10^{-4}\) for
every \(N\), 95th percentile \(2\times10^{-3}\) and maximum \(3\times10^{-2}\)
over the \(250\) chains; the two \(Q_1\) constructions on \(16^3\) cells have
\(5.7\)--\(6.8\%\). For \(f=0\) and end data in the port space the solution
contains only the port modes, so the interior port truncation is exact
and the Ms-LegONet error is that of the learned space alone. The gap is
discretization and, as for \(K^\ast\)
in \Cref{subsec:2d-diffusion}, it enters \(\mathbf D\) squared to leading
order: \(g^{T}(\mathbf D_h-\mathbf D)g=\|u_g-u_{g,h}\|_{a_\kappa}^2\) for a
conforming Galerkin solution \(u_{g,h}\) with the same end data \(g\).
The \(Q_1\) blocks impose the sine traces at the face nodes, so the
identity holds up to the interpolation error of the trace; the Frobenius
error of \(\mathbf D\) aggregates it over the \(32\) ports. The \(Q_1\)
field error is \(O(h)\), \(24\)--\(26\%\) at \(16^3\) and \(12\)--\(13\%\)
at \(32^3\), consistent with the \(6\%\) and \(1.6\%\) of
\Cref{tab:3d-eff}; the learned block is not bound by the \(Q_1\) floor of the \(16^3\) mesh, and the same square
turns its port-response energy error of about \(1\%\) into \(10^{-4}\).
Reaching \(10^{-4}\) with \(Q_1\) would take a field error of \(1\%\), by
extrapolation about \(400^3\) cells per primitive; one halving of \(h\) already costs \(25\) s
per coefficient at \(N=16\) and \(98.5\) s at \(N=64\), \(827\) and \(904\)
times Ms-LegONet, for \(1.6\%\) error.

\begin{table}[!htbp]
\centering
\caption{Many-query construction of the end-to-end response $\mathbf D$ ($32\times32$, $f=0$) of a chain of $N$ primitives: time per new coefficient and relative Frobenius error against the exact continuum $\mathbf D$, 50 test chains per $N$. $Q_1$ blocks: the Ms-LegONet construction with $Q_1$ local solves in place of the network, batched over the primitives. Monolithic $Q_1$: one algebraic multigrid hierarchy of the chain and 32 end-face solves by preconditioned CG to a relative residual of $10^{-8}$ (a tolerance of $10^{-4}$ lowers these times by about $40\%$). Speed-ups are computed from unrounded times.}
\label{tab:3d-eff}
\begin{adjustbox}{max width=\textwidth}
\begin{tabular}{lcccccc}
\toprule
 & \multicolumn{5}{c}{time per coefficient [s], $N$} & relative error of $\mathbf D$\\
\cmidrule(lr){2-6}
method & $4$ & $8$ & $16$ & $32$ & $64$ & per-$N$ medians, range\\
\midrule
\textbf{Ms-LegONet} & \textbf{0.016} & \textbf{0.020} & \textbf{0.031} & \textbf{0.053} & \textbf{0.109} & \textbf{$(6.7$--$9.9)\times10^{-5}$}\\
$Q_1$ blocks, $16^3$ & 0.508 & 0.858 & 1.604 & 3.186 & 6.352 & 0.057--0.068\\
monolithic $Q_1$, $16^3$ & 0.961 & 1.781 & 3.529 & 5.597 & 11.53 & 0.057--0.068\\
monolithic $Q_1$, $32^3$ & 6.34 & 13.8 & 25.3 & 47.3 & 98.5 & 0.015--0.018\\
\midrule
speed-up vs.\ $Q_1$ blocks, $16^3$ & 31 & 43 & 52 & 60 & 58 & \\
speed-up vs.\ monolithic $Q_1$, $16^3$ & 59 & 89 & 115 & 106 & 106 & \\
speed-up vs.\ monolithic $Q_1$, $32^3$ & 388 & 688 & 827 & 893 & 904 & \\
\bottomrule
\end{tabular}
\end{adjustbox}
\end{table}

\section{Conclusions}
\label{sec:conclusion}

Ms-LegONet is a plug-and-play multiscale block whose network amortizes
one local task, the coefficient-dependent harmonic space of a patch. One
network per dimension, trained once on \(1.3\times10^5\) and \(5{,}500\)
single-patch labels, serves, within its coefficient class, every layout,
source, boundary datum, layout size from \(2\) to \(4{,}096\) primitives
and composition depth. Everything downstream is deterministic linear algebra, and the
global solution is conforming Galerkin under the localized bound of
\Cref{prop:localized-bound}.

Three results carry the paper. First, the learned space is as accurate as
the local problem allows: in three dimensions the block stays within
\(0.1\) percentage points of the exact local spaces in median at every
chain length, and in two dimensions the root-sum-square of the local
best-approximation errors bounds the global error with effectivity within
\(10\%\) of one (\Cref{tab:2d-global}). Second, the local coarse
spaces of the classical baselines, not their meshes, limit them: at equal
unknowns the median energy error of Ms-LegONet is \(2.4\) to \(7.6\) times
below GMsFEM-32 \(+\,p_f\) and \(1.4\) to \(2.3\) times below NOEM in two
dimensions and about \(3\) times below the strongest classical variant in three
(\Cref{tab:2d-baselines-acc,tab:3d-baselines}), and an end-to-end FNO
trained on mosaics of the same patches reaches \(38\%\) \(L^2\) error. Third, the
cost gain is largest in three dimensions: the many-query port response of a
chain is computed \(31\) to \(115\) times faster than with finite element
blocks or a monolithic solve on the same mesh, and up to \(900\) times
faster than the monolithic solve on \(32^3\) cells, at \(10^{-4}\)
relative error against their \(6\%\) (\Cref{tab:3d-eff}); a monolithic
solve of that accuracy costs three to four orders of magnitude more. In
two dimensions the block is faster than a fine-grid direct solve while
reproducing the effective conductivity to within \(0.3\%\)
(\Cref{tab:2d-homog}).

The limitations we observe point to the next steps. In three dimensions
the remaining error is the port truncation of the \(32\) sine modes; a larger port
space removes it, and coefficients varying in all three directions will
require ports on all six faces, with the composition and the bound
unchanged. In two dimensions the learned space is least
accurate on a few high-contrast primitives, where lift and correction
nearly cancel; the localized bound identifies them, and with a computable
surrogate for \(\eta_{\mathrm{proj}}\), such as the residual of the learned
port responses, they could be solved exactly. The block is trained for a
coefficient class on a fixed primitive mesh, \(16\times16\) cells in two
dimensions and \(16^3\) in three; the test primitives of
\Cref{subsec:2d-diffusion,subsec:3d-layered} are new draws from that
class, while every layout, source and boundary datum is new. A class
outside it, a new microstructure family, a higher contrast or another
primitive resolution, is reached by fine-tuning the same block on
single-patch labels of the new class, at a data cost independent of
layout and domain size and with ports, composition and bound unchanged;
the bound marks the primitives on which the current block falls short.
Two extensions follow from the construction. A block that also predicts
the responses to a few fixed source modes would remove every local
factorization from the online stage; in a general medium, where \(p_f\)
is computed on the local mesh, it would also replace the discretization
error of \(p_f\), which then bounds the global accuracy
(\Cref{app:q1-pf}), by a learning error. And
because the learned object is a space, the block can serve as the local
component of other coarse-space methods, such as GMsFEM or localized
reduced-basis methods, and any symmetric coercive problem with a local
harmonic-extension structure, such as linear elasticity with vector-valued
traces, admits the same design.
\appendix

\section{Technical properties of the local block}
\label{app:local-block}

\subsection{Affine pullback and harmonic extensions}
\label{app:affine-pullback}

With the notation of \Cref{subsec:primitive}, let
\(F_i(\widehat{\mathbf x})=\mathbf B_i\widehat{\mathbf x}+\mathbf c_i\)
with \(\mathbf B_i\in\mathbb R^{d\times d}\) nonsingular. For
\(\widehat w=w\circ F_i\) and \(\widehat v=v\circ F_i\),
\begin{equation}
a_i^{\kappa_i}(w,v)
=
\int_{\widehat\omega}
\nabla_{\widehat{\mathbf x}}\widehat w^{\,T}
\widehat{\mathbf K}_i
\nabla_{\widehat{\mathbf x}}\widehat v
\,\mathrm d\widehat{\mathbf x},
\qquad
\widehat{\mathbf K}_i(\widehat{\mathbf x})
=
|\det\mathbf B_i|\,
\kappa_i(F_i(\widehat{\mathbf x}))
\mathbf B_i^{-1}\mathbf B_i^{-T}.
\label{eq:affine-pulled-coefficient}
\end{equation}
For translations, \(\mathbf B_i=\mathbf I\) and
\(\widehat{\mathbf K}_i=(\kappa_i\circ F_i)\mathbf I\).

\subsection{Energy orthonormalization}

\begin{proposition}[Subspace discrepancy]
\label{prop:app-energy-subspace-invariance}
Let \(\mathbf A_i(\kappa_i)\) be symmetric positive definite and
\(\widetilde{\mathbf Q}_{\theta,i}\in\mathbb R^{n_0\times q}\) of
full column rank. Then \(\mathbf Q_{\theta,i}\) of
\eqref{eq:energy-orthonormalized-neural-basis} satisfies
\begin{equation}
\mathbf Q_{\theta,i}^{T}
\mathbf A_i(\kappa_i)
\mathbf Q_{\theta,i}
=
\mathbf I_q,
\qquad
\operatorname{Range}(\mathbf Q_{\theta,i})
=
\operatorname{Range}(\widetilde{\mathbf Q}_{\theta,i}).
\label{eq:app-energy-orthonormality}
\end{equation}
If \(\mathbf Q_i^\star\) is any
\(\mathbf A_i(\kappa_i)\)-orthonormal basis of
\(\mathcal V_{i,q}^{\star}(\kappa_i)\) and
\(\{\vartheta_{r}\}_{r=1}^q\) are the principal angles between
\(\mathcal V_{i,q}^{\star}(\kappa_i)\) and
\(\mathcal V_{\theta,i}(\kappa_i)\) in the
\(\mathbf A_i(\kappa_i)\)-inner product, then
\begin{equation}
\cos\vartheta_{r}
=
\sigma_{r}
\left(
{\mathbf Q_i^\star}^{T}
\mathbf A_i(\kappa_i)
\mathbf Q_{\theta,i}
\right),
\qquad
r=1,\ldots,q,
\label{eq:app-energy-principal-angles}
\end{equation}
and
\begin{equation}
\sum_{r=1}^{q}\sin^2\vartheta_{r}
=
q-
\left\|
{\mathbf Q_i^\star}^{T}
\mathbf A_i(\kappa_i)
\mathbf Q_{\theta,i}
\right\|_F^2.
\label{eq:app-grassmann-discrepancy}
\end{equation}
The right side of \eqref{eq:app-grassmann-discrepancy} depends on the
two subspaces alone.
\end{proposition}

\begin{proof}
Suppressing \(i\) and \(\kappa_i\),
\(\mathbf G_\theta=\widetilde{\mathbf Q}_\theta^{T}\mathbf A\widetilde{\mathbf Q}_\theta\)
is symmetric positive definite since
\(\widetilde{\mathbf Q}_\theta\) has full column rank, so
\(\mathbf Q_\theta^{T}\mathbf A\mathbf Q_\theta=\mathbf G_\theta^{-1/2}\mathbf G_\theta\mathbf G_\theta^{-1/2}=\mathbf I_q\),
and the nonsingular factor
\(\mathbf G_\theta^{-1/2}\) preserves the column space.
Under the isometry \(\mathbf v\mapsto\mathbf A^{1/2}\mathbf v\)
from the \(\mathbf A\)-inner product to the Euclidean one,
\(\mathbf A^{1/2}\mathbf Q_i^\star\) and
\(\mathbf A^{1/2}\mathbf Q_\theta\) have Euclidean-orthonormal columns, and
the singular-value characterization of principal angles
\cite{bjorck1973numerical,li2006matrix,absil2008optimization}
gives \eqref{eq:app-energy-principal-angles}; summing
\(\sin^2\vartheta_{r}=1-\cos^2\vartheta_{r}\) gives
\eqref{eq:app-grassmann-discrepancy}.
\end{proof}

\subsection{Local Galerkin projection}

\begin{proposition}[Local Galerkin projection]
\label{prop:app-local-response-projection}
For any \(g\in\mathcal B_{\Gamma_i}\),
\eqref{eq:mslegonet-local-galerkin} admits a unique solution,
\begin{equation}
a_i^{\kappa_i}
\left(
\mathcal C_{i,h}^{\kappa_i}g-w_{\theta,i}(g),
v_\theta
\right)
=
0
\qquad
\forall v_\theta\in\mathcal V_{\theta,i}(\kappa_i),
\label{eq:app-local-galerkin-orthogonality}
\end{equation}
and \eqref{eq:local-response-best-approximation} holds. If
\(\mathcal R_{i,h}(\kappa_i)\subseteq\mathcal V_{\theta,i}(\kappa_i)\),
then
\(u_{\theta,i}(g)=\mathcal H_{i,h}^{\kappa_i}g\).
\end{proposition}

\begin{proof}
As \(a_i^{\kappa_i}\) is symmetric positive definite on
\(\mathcal V_{\theta,i}(\kappa_i)\subseteq V_{i,h}^0\), the reduced
problem is uniquely solvable. The exact extension satisfies
\(a_i^{\kappa_i}\left(E_i g+\mathcal C_{i,h}^{\kappa_i}g,v_\theta\right)=0\)
for all \(v_\theta\in\mathcal V_{\theta,i}(\kappa_i)\); subtracting
\eqref{eq:mslegonet-local-galerkin} gives
\eqref{eq:app-local-galerkin-orthogonality}. Thus
\(w_{\theta,i}(g)\) is the \(a_i^{\kappa_i}\)-orthogonal projection of
\(\mathcal C_{i,h}^{\kappa_i}g\) onto
\(\mathcal V_{\theta,i}(\kappa_i)\), and
\eqref{eq:local-response-best-approximation} follows. Since
\(\mathcal C_{i,h}^{\kappa_i}g\in\mathcal R_{i,h}(\kappa_i)\), the
inclusion makes the projection exact.
\end{proof}

\section{Recursive composition}
\label{app:recursive-composition-algebra}

Notation follows \Cref{subsec:composition,subsec:source-global}.

\subsection{Parent composition}

\begin{proposition}[Parent composition]
\label{prop:app-parent-composition}
Assume \eqref{eq:kappa-positive}, connectedness of \(\omega_P\),
conforming identification of the child port spaces on shared
interfaces, and injectivity of the parametrization of
\(\Lambda_P^{\mathrm{int}}\). Then, for every
\(g_P\in\mathcal B_{\Gamma_P}\),
\eqref{eq:parent-energy-minimization} has a unique solution
\(\lambda_P^\star(g_P)\). The parent field
\(\mathcal H_{\theta,P}g_P\) of
\eqref{eq:parent-response-map} lies in \(H^1(\omega_P)\), has trace
\(g_P\) on \(\Gamma_P\), satisfies
\eqref{eq:interface-equilibrium}, and depends linearly on \(g_P\).
\end{proposition}

\begin{proof}
By linearity of the child response maps, \(J_P(\cdot;g_P)\) is
quadratic on \(\Lambda_P^{\mathrm{int}}\), with quadratic form in the direction
\(\mu\in\Lambda_P^{\mathrm{int}}\)
\begin{equation}
\sum_{D\in\mathcal D(P)}
a_D^{\kappae}
\left(
\mathcal H_{\theta,D}
\mathcal M_{P\to D}^{I}\mu,
\mathcal H_{\theta,D}
\mathcal M_{P\to D}^{I}\mu
\right).
\label{eq:app-parent-energy-hessian}
\end{equation}
If this nonnegative quantity vanishes, uniform positivity of
\(\kappae\) makes every induced child field constant on its connected
child block; trace conformity makes these constants agree across adjacent
children, and the zero trace on \(\Gamma_P\) with connectedness of
\(\omega_P\) makes the assembled field vanish. Injectivity of the skeleton
parametrization then gives \(\mu=0\): the quadratic form is positive
definite and the minimizer unique.
Matching child traces on internal interfaces put the parent field in
\(H^1(\omega_P)\) \cite{brenner2008mathematical}, with trace \(g_P\) on
\(\Gamma_P\) by construction. Differentiating \(J_P\) at
the minimizer gives \eqref{eq:interface-equilibrium}, linear in
\(g_P\); uniqueness makes the parent map linear.
\end{proof}

\subsection{Static condensation}

Choose bases for \(\mathcal B_{\Gamma_P}\) and \(\Lambda_P^{\mathrm{int}}\)
with coefficient vectors \(\boldsymbol\eta\) and \(\boldsymbol\xi\), let
\(\mathbf M_D^E\) and \(\mathbf M_D^I\) be the matrices of the trace maps
in \eqref{eq:child-induced-trace}, so that the child trace coefficients
are \(\mathbf c_D=\mathbf M_D^E\boldsymbol\eta+\mathbf M_D^I\boldsymbol\xi\),
and let \(\mathbf S_{\theta,D}\) be the matrix of
\(s_{\theta,D}^{\kappae}\) in the child port basis. Then
\begin{equation}
J_P
=
\frac12
\begin{bmatrix}
\boldsymbol\eta\\
\boldsymbol\xi
\end{bmatrix}^{T}
\begin{bmatrix}
\mathbf K_{EE} & \mathbf K_{EI}\\
\mathbf K_{EI}^{T} & \mathbf K_{II}
\end{bmatrix}
\begin{bmatrix}
\boldsymbol\eta\\
\boldsymbol\xi
\end{bmatrix},
\qquad
\mathbf K_{XY}
=
\sum_{D\in\mathcal D(P)}
(\mathbf M_D^X)^T
\mathbf S_{\theta,D}
\mathbf M_D^Y,
\label{eq:app-parent-block-energy}
\end{equation}
for \(X,Y\in\{E,I\}\).

\begin{proposition}[Parent Schur complement]
\label{prop:app-parent-schur}
Under the assumptions of \Cref{prop:app-parent-composition},
\(\mathbf K_{II}\) is symmetric positive definite,
\begin{equation}
\boldsymbol\xi^\star
=
-\mathbf K_{II}^{-1}\mathbf K_{EI}^{T}\boldsymbol\eta,
\qquad
\mathbf S_{\theta,P}
=
\mathbf K_{EE}
-
\mathbf K_{EI}
\mathbf K_{II}^{-1}
\mathbf K_{EI}^{T},
\label{eq:app-parent-schur-complement}
\end{equation}
and
\(J_P\bigl(\lambda_P^\star;g_P\bigr)=\frac12\boldsymbol\eta^T\mathbf S_{\theta,P}\boldsymbol\eta\).
\end{proposition}

\begin{proof}
\(\mathbf K_{II}\) is the matrix of the positive-definite form
\eqref{eq:app-parent-energy-hessian}. Differentiating
\eqref{eq:app-parent-block-energy} with respect to
\(\boldsymbol\xi\) gives
\(\mathbf K_{II}\boldsymbol\xi^\star=-\mathbf K_{EI}^{T}\boldsymbol\eta\);
substituting into \eqref{eq:app-parent-block-energy}
gives \eqref{eq:app-parent-schur-complement}.
\end{proof}

\subsection{Particular--harmonic decomposition}

\begin{proposition}[Particular--harmonic decomposition]
\label{prop:app-particular-harmonic-decomposition}
Let \(u_h\in\overline V_h^{g_\Omega}\) solve
\eqref{eq:global-source-discrete-weak-form} and \(p_f\) be given
by \eqref{eq:local-particular-problem}. Then
\(r_h:=u_h-p_f\) has trace \(g_\Omega\) on \(\partial\Omega\) and is
discretely \(\kappa_i\)-harmonic on every primitive:
\begin{equation}
a_i^{\kappa_i}
\left(
r_h|_{\omega_i},v_i
\right)
=
0
\qquad
\forall v_i\in V_{i,h}^0,
\qquad
i=1,\ldots,N_H.
\label{eq:app-remainder-local-harmonicity}
\end{equation}
\end{proposition}

\begin{proof}
For \(v_i\in V_{i,h}^0\), the zero extension \(\mathcal E_i^0v_i\)
lies in \(V_h\); testing
\eqref{eq:global-source-discrete-weak-form} with it and using
\eqref{eq:local-particular-problem} gives
\(a_i^{\kappa_i}\left(u_h|_{\omega_i},v_i\right)=\ell_f(\mathcal E_i^0v_i)=\ell_{f,i}(v_i)=a_i^{\kappa_i}\bigl(p_i(f),v_i\bigr)\).
The zero extensions of \(p_j(f)\), \(j\ne i\), vanish in \(\omega_i\),
so subtracting gives
\eqref{eq:app-remainder-local-harmonicity}. Since
\(p_f\in V_h\subset H_0^1(\Omega)\), subtracting \(p_f\) preserves
the trace.
\end{proof}
\section{Proof of the localized bound}
\label{app:localized-bound-proof}

\begin{proof}[Proof of \Cref{prop:localized-bound}]
By \eqref{eq:local-response-best-approximation}, the learned
correction of a trace \(g\) is the projection of the exact correction,
\(w_{\theta,i}(g)=\Pi_{\theta,i}\mathcal C_{i,h}^{\kappa_i}g\).
Define \(v\) on every primitive by
\(v|_{\omega_i}:=p_i(f)+E_ig_i^{\star}+\Pi_{\theta,i}w_i^{\star}\).
Since the skeleton traces \(g_i^{\star}\) are admissible and
conforming, this function belongs to the affine Ms-LegONet trial
space \(p_f+V_{\mathrm{ms},\theta}^{0}\) (shifted by the boundary
response if \(g_\Omega\neq0\)). Since \(u_{\mathrm{ms},\theta}\) is
the \(a_{\kappae}\)-orthogonal projection of \(u_h\) onto this affine
space,
\(\|u_h-u_{\mathrm{ms},\theta}\|_{a_{\kappae}}\le\|u_h-v\|_{a_{\kappae}}\le\|u_h-u_{\mathrm{ms}}^{\star}\|_{a_{\kappae}}+\|u_{\mathrm{ms}}^{\star}-v\|_{a_{\kappae}}\).
On \(\omega_i\), \(u_{\mathrm{ms}}^{\star}-v\) equals
\((I-\Pi_{\theta,i})w_i^{\star}\), so the second term equals
\(\bigl(\sum_i\|(I-\Pi_{\theta,i})w_i^{\star}\|_{a_i^{\kappa_i}}^2\bigr)^{1/2}\)\({}=\eta_{\mathrm{proj}}\), which gives the first inequality of \eqref{eq:localized-bound}.
The second inequality follows from
\(\|(I-\Pi_{\theta,i})w\|_{a_i^{\kappa_i}}\le\sin\vartheta_{\max,i}\|w\|_{a_i^{\kappa_i}}\)
for every \(w\in\mathcal V_{i,q}^{\star}(\kappa_i)\). This inequality
is the characterization of the largest principal angle between two
subspaces of equal dimension \cite{bjorck1973numerical}. In the
full-trace case the exact-block construction is the exact static
condensation of \eqref{eq:global-source-discrete-weak-form}
(\Cref{app:recursive-composition-algebra}), so
\(u_{\mathrm{ms}}^{\star}=u_h\).
\end{proof}

\section{Training settings}
\label{app:training-settings}

The input is normalized with \(s_\kappa=\log10^4\) in both dimensions.
The two-dimensional network is a residual convolutional encoder with four
resolution levels and a pointwise MLP decoder that returns the \(q=63\)
nodal zero-trace fields. The three-dimensional network is a residual
convolutional encoder followed by directional pooling and three MLP heads
that return the one-dimensional factors of the separated representation
of \Cref{subsec:3d-layered}.

The loss uses $q_{\mathrm{tail}}=8$, $\rho=0.2$ and $M=24$ probe traces per patch in two dimensions and $M=12$ in three. The weights start at $\beta=(1,1,0)$; early in training $\beta_{\mathrm{space}}$ falls linearly to $0.5$ and $\beta_{\mathrm{phys}}$ rises linearly to $1$, after which they stay fixed. Optimization uses Adam with gradient-norm clipping at $1$ and minibatches balanced over families and contrasts, with batch $96$, peak learning rate $2\times10^{-4}$ and $600$k updates in two dimensions and batch $32$, peak learning rate $3\times10^{-4}$ and $120$k updates in three. All remaining settings are in the public repository.

\section{Global error with the particular fields on the \texorpdfstring{$Q_1$}{Q1} mesh}
\label{app:q1-pf}

In \Cref{subsec:3d-layered} the local particular fields \(p_i(f)\) are
computed with the exact modal solver. \Cref{tab:q1-pf} reports the global
energy error of Ms-LegONet on the \(200\) test chains with \(N=16\) and
\(f=1\) when \(p_i(f)\) is instead obtained by a \(Q_1\) solve on
\(16^3\), \(32^3\) or \(64^3\) cells per primitive, with the learned
harmonic space and the composition unchanged. The error follows the
\(Q_1\) error of \(p_f\), first order in \(h\) and equal to that of the
monolithic \(Q_1\) solution on the same mesh, plus the port floor of
\(4.4\%\) approximately in quadrature (with an approximate \(p_f\) the two
errors are no longer exactly orthogonal): \(\sqrt{7.27^2+4.41^2}=8.5\)
against the measured \(8.76\%\) on \(64^3\) cells. On the \(16^3\) mesh the median error by
contrast group is \(21.2\%\) at \(\chi=10^2\), \(31.4\%\) at \(10^3\),
\(42.1\%\) at \(10^4\) and \(27.7\%\) on the mixed chains. Since
\(p_i(f)\) is local and free of learned quantities, its mesh can be chosen
independently of the network's \(16^3\) input mesh.

\begin{table}[!htbp]
\centering
\caption{Relative energy error (\%) of Ms-LegONet against the continuum reference on the 200 test chains with $N=16$, $f=1$, when the particular fields $p_i(f)$ are computed by $Q_1$ solves on the stated local mesh instead of the modal solver; the monolithic $Q_1$ rows of \Cref{tab:3d-baselines} are repeated for comparison.}
\label{tab:q1-pf}
\begin{tabular}{lrr}
\toprule
method & P50 & P95\\
\midrule
Ms-LegONet, exact $p_f$ (\Cref{tab:3d-baselines}) & 4.51 & 7.10\\
Ms-LegONet, $Q_1$ $p_f$ on $16^3$ & 29.59 & 45.99\\
Ms-LegONet, $Q_1$ $p_f$ on $32^3$ & 15.44 & 23.81\\
Ms-LegONet, $Q_1$ $p_f$ on $64^3$ & 8.76 & 12.97\\
\addlinespace
monolithic $Q_1$, $16^3$ & 28.94 & 45.78\\
monolithic $Q_1$, $32^3$ & 14.52 & 22.97\\
monolithic $Q_1$, $64^3$ & 7.27 & 11.50\\
\bottomrule
\end{tabular}
\end{table}

\section*{Data and code availability}
The implementation, the training-data generators, the pretrained networks and the scripts that reproduce every table and figure are available at \url{https://github.com/Yooki-YueqiWang/Ms-LegONet}.

\section*{Acknowledgments}
During the preparation of this work the authors used Claude (Anthropic) to improve the readability and language of the manuscript. After using this tool, the authors reviewed and edited the content as needed. The authors assume responsibility for all content.

We gratefully acknowledge the support from DARPA MATHBAC program under grant number: HR0011-26-3-E152, the National Science Foundation (DMS-2533878, DMS-2053746, DMS-2134209, ECCS-2328241, CBET-2347401 and OAC-2311848), and U.S.~Department of Energy (DOE) Office of Science Advanced Scientific Computing Research program under the "Uncertainty Quantification for Multifidelity Operator Learning (MOLUcQ)" project (Project No. 81739), DE-SC0023161, the SciDAC LEADS Institute, and DOE–Fusion Energy Science, under grant number: DE-SC0024583.

\bibliographystyle{siamplain}
\bibliography{reference2}

@article{hou1997msfem,
  author         = {Hou, Thomas Y. and Wu, Xiao-Hui},
  title          = {A Multiscale Finite Element Method for Elliptic Problems in Composite Materials and Porous Media},
  journal        = {Journal of Computational Physics},
  volume         = {134},
  number         = {1},
  pages          = {169--189},
  year           = {1997},
  doi            = {10.1006/jcph.1997.5682}
}

@article{efendiev2011spectral,
  author         = {Efendiev, Yalchin and Galvis, Juan and Wu, Xiao-Hui},
  title          = {Multiscale finite element methods for high-contrast problems using local spectral basis functions},
  journal        = {Journal of Computational Physics},
  volume         = {230},
  number         = {4},
  pages          = {937--955},
  year           = {2011},
  doi            = {10.1016/j.jcp.2010.09.026}
}

@article{efendiev2013gmsfem,
  author         = {Efendiev, Yalchin and Galvis, Juan and Hou, Thomas Y.},
  title          = {Generalized multiscale finite element methods ({GMsFEM})},
  journal        = {Journal of Computational Physics},
  volume         = {251},
  pages          = {116--135},
  year           = {2013},
  doi            = {10.1016/j.jcp.2013.04.045}
}

@article{abdulle2012rbfehmm,
  author         = {Abdulle, Assyr and Bai, Yun},
  title          = {Reduced basis finite element heterogeneous multiscale method for high-order discretizations of elliptic homogenization problems},
  journal        = {Journal of Computational Physics},
  volume         = {231},
  number         = {21},
  pages          = {7014--7036},
  year           = {2012},
  doi            = {10.1016/j.jcp.2012.02.019}
}

@article{ohlberger2015lrbms,
  author         = {Ohlberger, Mario and Schindler, Felix},
  title          = {Error Control for the Localized Reduced Basis Multiscale Method with Adaptive On-Line Enrichment},
  journal        = {SIAM Journal on Scientific Computing},
  volume         = {37},
  number         = {6},
  pages          = {A2865--A2895},
  year           = {2015},
  doi            = {10.1137/151003660}
}

@article{buhr2017arbilomod,
  author         = {Buhr, Andreas and Engwer, Christian and Ohlberger, Mario and Rave, Stephan},
  title          = {{ArbiLoMod}, a Simulation Technique Designed for Arbitrary Local Modifications},
  journal        = {SIAM Journal on Scientific Computing},
  volume         = {39},
  number         = {4},
  pages          = {A1435--A1465},
  year           = {2017},
  doi            = {10.1137/15M1054213}
}

@article{lu2021deeponet,
  author         = {Lu, Lu and Jin, Pengzhan and Pang, Guofei and Zhang, Zhongqiang and Karniadakis, George Em},
  title          = {Learning nonlinear operators via {DeepONet} based on the universal approximation theorem of operators},
  journal        = {Nature Machine Intelligence},
  volume         = {3},
  number         = {3},
  pages          = {218--229},
  year           = {2021},
  doi            = {10.1038/s42256-021-00302-5},
  url            = {https://www.nature.com/articles/s42256-021-00302-5}
}

@inproceedings{li2021fno,
  author         = {Li, Zongyi and Kovachki, Nikola B. and Azizzadenesheli, Kamyar and Liu, Burigede and Bhattacharya, Kaushik and Stuart, Andrew and Anandkumar, Anima},
  title          = {Fourier Neural Operator for Parametric Partial Differential Equations},
  booktitle      = {International Conference on Learning Representations},
  year           = {2021},
  url            = {https://openreview.net/forum?id=c8P9NQVtmnO}
}

@article{kovachki2023neuraloperator,
  author         = {Kovachki, Nikola and Li, Zongyi and Liu, Burigede and Azizzadenesheli, Kamyar and Bhattacharya, Kaushik and Stuart, Andrew and Anandkumar, Anima},
  title          = {Neural Operator: Learning Maps Between Function Spaces With Applications to {PDEs}},
  journal        = {Journal of Machine Learning Research},
  volume         = {24},
  number         = {89},
  pages          = {1--97},
  year           = {2023},
  url            = {https://www.jmlr.org/papers/v24/21-1524.html}
}

@article{li2024pino,
  author         = {Li, Zongyi and Zheng, Hongkai and Kovachki, Nikola and Jin, David and Chen, Haoxuan and Liu, Burigede and Azizzadenesheli, Kamyar and Anandkumar, Anima},
  title          = {Physics-Informed Neural Operator for Learning Partial Differential Equations},
  journal        = {ACM/IMS Journal of Data Science},
  volume         = {1},
  number         = {3},
  pages          = {1--27},
  year           = {2024},
  doi            = {10.1145/3648506},
  url            = {https://doi.org/10.1145/3648506}
}

@misc{lutjens2022mno,
  author        = {L{\"u}tjens, Bj{\"o}rn and Crawford, Catherine H. and Watson, Campbell D. and Hill, Christopher and Newman, Dava},
  title         = {Multiscale Neural Operator: Learning Fast and Grid-independent {PDE} Solvers},
  year          = {2022},
  eprint        = {2207.11417},
  archiveprefix = {arXiv},
  primaryclass  = {cs.LG}
}

@inproceedings{he2024mgno,
  author         = {He, Juncai and Liu, Xinliang and Xu, Jinchao},
  title          = {{MgNO}: Efficient Parameterization of Linear Operators via Multigrid},
  booktitle      = {International Conference on Learning Representations},
  year           = {2024},
  url            = {https://proceedings.iclr.cc/paper_files/paper/2024/hash/eb3c8135137c8a60425a0320869ad87e-Abstract-Conference.html}
}

@article{li2024lno,
  author         = {Li, Hongyu and Ye, Ximeng and Jiang, Peng and Qin, Guoliang and Wang, Tiejun},
  title          = {Local neural operator for solving transient partial differential equations on varied domains},
  journal        = {Computer Methods in Applied Mechanics and Engineering},
  volume         = {427},
  pages          = {117062},
  year           = {2024},
  doi            = {10.1016/j.cma.2024.117062},
  url            = {https://doi.org/10.1016/j.cma.2024.117062}
}

@article{fang2024boundary,
  author         = {Fang, Zhiwei and Wang, Sifan and Perdikaris, Paris},
  title          = {Learning Only on Boundaries: A Physics-Informed Neural Operator for Solving Parametric Partial Differential Equations in Complex Geometries},
  journal        = {Neural Computation},
  volume         = {36},
  number         = {3},
  pages          = {475--498},
  year           = {2024},
  doi            = {10.1162/neco_a_01647},
  url            = {https://doi.org/10.1162/neco_a_01647}
}

@article{chen2025locality,
  author         = {Chen, Jiangce and Xu, Wenzhuo and Xu, Zeda and Grande Guti{\'e}rrez, Noelia and Narra, Sneha Prabha and McComb, Christopher},
  title          = {Enforcing the principle of locality for physical simulations with neural operators},
  journal        = {Journal of Computational Physics},
  volume         = {538},
  pages          = {114131},
  year           = {2025},
  doi            = {10.1016/j.jcp.2025.114131},
  url            = {https://doi.org/10.1016/j.jcp.2025.114131}
}

@inproceedings{huang2026sni,
  author         = {Huang, Jianing and Zhang, Kaixuan and Wu, Youjia and Cheng, Ze},
  title          = {Operator Learning with Domain Decomposition for Geometry Generalization in {PDE} Solving},
  booktitle      = {International Conference on Learning Representations},
  year           = {2026},
  url            = {https://proceedings.iclr.cc/paper_files/paper/2026/hash/710445227fa8c1b6a9ceada902dd4741-Abstract-Conference.html}
}

@article{he2026selfcompose,
  author         = {He, Juncai and Liu, Xinliang and Xu, Jinchao},
  title          = {Self-composing neural operators for high-frequency and multiscale {PDE} surrogates},
  journal        = {Journal of Computational Physics},
  volume         = {565},
  pages          = {115189},
  year           = {2026},
  doi            = {10.1016/j.jcp.2026.115189},
  url            = {https://doi.org/10.1016/j.jcp.2026.115189}
}

@article{ouyang2026noem,
  author         = {Ouyang, Weihang and Shin, Yeonjong and Liu, Si-Wei and Lu, Lu},
  title          = {{NOEM}: efficient and scalable finite element method enabled by reusable neural operators},
  journal        = {Nature Computational Science},
  volume         = {6},
  number         = {4},
  pages          = {417--429},
  year           = {2026},
  doi            = {10.1038/s43588-026-00974-2},
  url            = {https://www.nature.com/articles/s43588-026-00974-2}
}

@misc{zhang2026legonet,
  author        = {Zhang, Jiahao and Wang, Yueqi and Lin, Guang},
  title         = {{LegONet}: Plug-and-Play Structure-Preserving Neural Operator Blocks for Compositional {PDE} Learning},
  year          = {2026},
  eprint        = {2603.07882},
  archiveprefix = {arXiv},
  primaryclass  = {math.NA}
}

@book{efendiev2009multiscale,
  author         = {Efendiev, Yalchin and Hou, Thomas Y.},
  title          = {Multiscale Finite Element Methods: Theory and Applications},
  series         = {Surveys and Tutorials in the Applied Mathematical Sciences},
  volume         = {4},
  year           = {2009},
  publisher      = {Springer},
  address        = {New York, NY},
  doi            = {10.1007/978-0-387-09496-0},
  isbn           = {978-0-387-09495-3},
  url            = {https://link.springer.com/book/10.1007/978-0-387-09496-0}
}

@book{chung2023multiscale,
  author         = {Chung, Eric and Efendiev, Yalchin and Hou, Thomas Y.},
  title          = {Multiscale Model Reduction: Multiscale Finite Element Methods and Their Generalizations},
  series         = {Applied Mathematical Sciences},
  volume         = {212},
  year           = {2023},
  publisher      = {Springer},
  address        = {Cham},
  doi            = {10.1007/978-3-031-20409-8},
  isbn           = {978-3-031-20408-1},
  url            = {https://link.springer.com/book/10.1007/978-3-031-20409-8}
}

@book{evans2022partial,
  author         = {Evans, Lawrence C.},
  title          = {Partial Differential Equations},
  edition        = {Second},
  series         = {Graduate Studies in Mathematics},
  volume         = {19},
  year           = {2022},
  publisher      = {American Mathematical Society},
  address        = {Providence, RI},
  isbn           = {978-1-4704-6942-9},
  url            = {https://katalog.bibliothek.kit.edu/bib/1337996}
}

@book{brenner2008mathematical,
  author         = {Brenner, Susanne C. and Scott, L. Ridgway},
  title          = {The Mathematical Theory of Finite Element Methods},
  edition        = {Third},
  series         = {Texts in Applied Mathematics},
  volume         = {15},
  year           = {2008},
  publisher      = {Springer},
  address        = {New York, NY},
  doi            = {10.1007/978-0-387-75934-0},
  isbn           = {978-0-387-75933-3},
  url            = {https://link.springer.com/book/10.1007/978-0-387-75934-0}
}

@book{ern2004theory,
  author         = {Ern, Alexandre and Guermond, Jean-Luc},
  title          = {Theory and Practice of Finite Elements},
  series         = {Applied Mathematical Sciences},
  volume         = {159},
  year           = {2004},
  publisher      = {Springer},
  address        = {New York, NY},
  doi            = {10.1007/978-1-4757-4355-5},
  isbn           = {978-0-387-20574-8},
  url            = {https://link.springer.com/book/10.1007/978-1-4757-4355-5}
}

@article{bjorck1973numerical,
  author         = {Bj{\"o}rck, {\AA}ke and Golub, Gene H.},
  title          = {Numerical Methods for Computing Angles Between Linear Subspaces},
  journal        = {Mathematics of Computation},
  volume         = {27},
  number         = {123},
  pages          = {579--594},
  year           = {1973},
  doi            = {10.1090/S0025-5718-1973-0348991-3}
}

@incollection{li2006matrix,
  author         = {Li, Ren-Cang},
  editor         = {Hogben, Leslie},
  title          = {Matrix Perturbation Theory},
  booktitle      = {Handbook of Linear Algebra},
  chapter        = {15},
  pages          = {15-1--15-17},
  year           = {2006},
  publisher      = {Chapman and Hall/CRC},
  address        = {Boca Raton, FL},
  isbn           = {978-1-58488-510-8},
  url            = {https://bpb-us-e1.wpmucdn.com/websites.uta.edu/dist/7/5059/files/2021/06/li2007HLA.pdf}
}

@book{absil2008optimization,
  author         = {Absil, P.-A. and Mahony, Robert and Sepulchre, Rodolphe},
  title          = {Optimization Algorithms on Matrix Manifolds},
  year           = {2008},
  publisher      = {Princeton University Press},
  address        = {Princeton, NJ},
  doi            = {10.1515/9781400830244},
  isbn           = {978-0-691-13298-3},
  url            = {https://press.princeton.edu/books/hardcover/9780691132983/optimization-algorithms-on-matrix-manifolds}
}

@book{quarteroni1999domain,
  author         = {Quarteroni, Alfio and Valli, Alberto},
  title          = {Domain Decomposition Methods for Partial Differential Equations},
  year           = {1999},
  publisher      = {Oxford University Press},
  address        = {Oxford},
  doi            = {10.1093/oso/9780198501787.001.0001},
  isbn           = {978-0-19-850178-7},
  url            = {https://academic.oup.com/book/53865}
}

@book{toselli2004domain,
  author         = {Toselli, Andrea and Widlund, Olof B.},
  title          = {Domain Decomposition Methods---Algorithms and Theory},
  series         = {Springer Series in Computational Mathematics},
  volume         = {34},
  year           = {2004},
  publisher      = {Springer},
  address        = {Berlin, Heidelberg},
  doi            = {10.1007/b137868},
  isbn           = {978-3-540-20696-5},
  url            = {https://link.springer.com/book/10.1007/b137868}
}

@article{rockafellar2000optimization,
  author         = {Rockafellar, R. Tyrrell and Uryasev, Stanislav},
  title          = {Optimization of Conditional Value-at-Risk},
  journal        = {Journal of Risk},
  volume         = {2},
  number         = {3},
  pages          = {21--41},
  year           = {2000},
  doi            = {10.21314/JOR.2000.038}
}

@misc{kingma2014adam,
  author        = {Kingma, Diederik P. and Ba, Jimmy},
  title         = {{Adam}: A Method for Stochastic Optimization},
  year          = {2014},
  eprint        = {1412.6980},
  archiveprefix = {arXiv},
  primaryclass  = {cs.LG}
}

@techreport{coons1967surfaces,
  author         = {Coons, Steven A.},
  title          = {Surfaces for Computer-Aided Design of Space Forms},
  number         = {MAC-TR-041},
  year           = {1967},
  month          = {June},
  institution    = {Massachusetts Institute of Technology, Project MAC},
  address        = {Cambridge, MA},
  url            = {https://hdl.handle.net/1721.1/149362},
  note           = {Also cataloged as MIT-LCS-TR-041}
}

@article{gordon1973transfinite,
  author         = {Gordon, William J. and Hall, Charles A.},
  title          = {Transfinite element methods: Blending-function interpolation over arbitrary curved element domains},
  journal        = {Numerische Mathematik},
  volume         = {21},
  number         = {2},
  pages          = {109--129},
  year           = {1973},
  doi            = {10.1007/BF01436298}
}

@book{bear2013dynamics,
  author         = {Bear, Jacob},
  title          = {Dynamics of Fluids in Porous Media},
  year           = {2013},
  publisher      = {Dover Publications},
  isbn           = {978-0-486-13180-1},
  url            = {https://search.worldcat.org/title/1136359830},
  note           = {Electronic reissue of the original 1972 work}
}

@book{torquato2002random,
  author         = {Torquato, Salvatore},
  title          = {Random Heterogeneous Materials: Microstructure and Macroscopic Properties},
  series         = {Interdisciplinary Applied Mathematics},
  volume         = {16},
  year           = {2002},
  publisher      = {Springer},
  address        = {New York, NY},
  doi            = {10.1007/978-1-4757-6355-3},
  isbn           = {978-0-387-95167-6},
  url            = {https://link.springer.com/book/10.1007/978-1-4757-6355-3}
}

@book{bensoussan2011asymptotic,
  author         = {Bensoussan, Alain and Lions, Jacques-Louis and Papanicolaou, George},
  title          = {Asymptotic Analysis for Periodic Structures},
  series         = {AMS Chelsea Publishing},
  volume         = {374},
  year           = {2011},
  publisher      = {American Mathematical Society},
  address        = {Providence, RI},
  isbn           = {978-0-8218-5324-5},
  url            = {https://www.ams.org/bookpages/chel-374},
  note           = {Corrected reprint of the 1978 original}
}

@book{jikov2012homogenization,
  author         = {Jikov, V. V. and Kozlov, S. M. and Oleinik, O. A.},
  title          = {Homogenization of Differential Operators and Integral Functionals},
  year           = {2012},
  publisher      = {Springer},
  address        = {Berlin, Heidelberg},
  doi            = {10.1007/978-3-642-84659-5},
  isbn           = {978-3-642-84659-5},
  url            = {https://link.springer.com/book/10.1007/978-3-642-84659-5},
  note           = {Electronic republication of the 1994 edition}
}

@book{ghanem2003stochastic,
  author         = {Ghanem, Roger G. and Spanos, Pol D.},
  title          = {Stochastic Finite Elements: A Spectral Approach},
  edition        = {Revised},
  year           = {2003},
  publisher      = {Dover Publications},
  address        = {Mineola, NY},
  isbn           = {978-0-486-42818-5},
  url            = {https://store.doverpublications.com/products/9780486428185}
}

@article{kuo2016application,
  author         = {Kuo, Frances Y. and Nuyens, Dirk},
  title          = {Application of Quasi-{Monte Carlo} Methods to Elliptic {PDEs} with Random Diffusion Coefficients: A Survey of Analysis and Implementation},
  journal        = {Foundations of Computational Mathematics},
  volume         = {16},
  number         = {6},
  pages          = {1631--1696},
  year           = {2016},
  doi            = {10.1007/s10208-016-9329-5}
}

@article{benner2015survey,
  author         = {Benner, Peter and Gugercin, Serkan and Willcox, Karen},
  title          = {A Survey of Projection-Based Model Reduction Methods for Parametric Dynamical Systems},
  journal        = {SIAM Review},
  volume         = {57},
  number         = {4},
  pages          = {483--531},
  year           = {2015},
  doi            = {10.1137/130932715}
}

@misc{rudikov2025gmsfemno,
  author        = {Rudikov, Alexander and Fanaskov, Vladimir and Stepanov, Sergei and Shan, Buzheng and Muravleva, Ekaterina and Efendiev, Yalchin and Oseledets, Ivan},
  title         = {Locally Subspace-Informed Neural Operators for Efficient Multiscale {PDE} Solving},
  year          = {2025},
  eprint        = {2505.16030},
  archiveprefix = {arXiv},
  primaryclass  = {math.NA}
}

@article{hou1999convergence,
  author         = {Hou, Thomas Y. and Wu, Xiao-Hui and Cai, Zhiqiang},
  title          = {Convergence of a multiscale finite element method for elliptic problems with rapidly oscillating coefficients},
  journal        = {Mathematics of Computation},
  volume         = {68},
  number         = {227},
  pages          = {913--943},
  year           = {1999},
  doi            = {10.1090/S0025-5718-99-01077-7}
}

@misc{zhang2026geometry,
  author        = {Zhang, Jiahao and Wang, Yueqi and Lin, Guang},
  title         = {Geometry-aware {LegONet} for {PDE} Learning on Arbitrary Domains},
  year          = {2026},
  eprint        = {2607.23069},
  archiveprefix = {arXiv}
}

@article{huynh2013scrbe,
  author         = {Huynh, Dinh Bao Phuong and Knezevic, David J. and Patera, Anthony T.},
  title          = {A static condensation reduced basis element method: approximation and a posteriori error estimation},
  journal        = {ESAIM: Mathematical Modelling and Numerical Analysis},
  volume         = {47},
  number         = {1},
  pages          = {213--251},
  year           = {2013},
  doi            = {10.1051/m2an/2012022}
}

@article{smetana2016optimal,
  author         = {Smetana, Kathrin and Patera, Anthony T.},
  title          = {Optimal local approximation spaces for component-based static condensation procedures},
  journal        = {SIAM Journal on Scientific Computing},
  volume         = {38},
  number         = {5},
  pages          = {A3318--A3356},
  year           = {2016},
  doi            = {10.1137/15M1009603}
}

\end{document}